\documentclass[11pt,a4paper]{article}

\usepackage[english]{babel}
\usepackage[utf8]{inputenc}
\usepackage{amsmath}
\usepackage{amsfonts}
\usepackage{graphicx}
\usepackage{epstopdf}
\usepackage[font=footnotesize, labelfont=bf]{caption}
\usepackage[figbotcap, loose]{subfigure}

\usepackage{siunitx}
\usepackage{tabularx}
\usepackage{booktabs}
\usepackage{multirow}
\usepackage{cleveref}

\usepackage{chngcntr}
\counterwithin*{equation}{section}

\usepackage{color}
\usepackage{soul}

\usepackage[square,numbers,sort&compress]{natbib}

\def\piergab #1{{\color{black}#1}}
\def\pier #1{{\color{red}#1}}
\def\pier #1{#1}
\newcommand{\EEE}{\color{black}}
\newcommand{\BBB}{\color{black}}
\def\pcol #1{{\color{red}#1}}
\def\pcol #1{#1}

\def\vec#1{\boldsymbol{#1}}

\begin{document}

\begin{center}
\textbf{\Large Mathematical analysis and simulation study of a phase-field model of prostate cancer growth with chemotherapy and antiangiogenic therapy effects}

\vspace{5mm}

{Pierluigi Colli}

{\footnotesize
{Dipartimento di Matematica, Universit\`a degli Studi di Pavia and IMATI-C.N.R., \\
Via Ferrata~5, 27100 Pavia, Italy\\
%and\\
%Istituto di Matematica Applicata e Tecnologie Informatiche ``Enrico Magenes'', CNR,\\
%Via Ferrata~1, 27100 Pavia, Italy\\
pierluigi.colli@unipv.it
}}

\vspace{2mm}

{Hector Gomez}

{\footnotesize
{School of Mechanical Engineering, Purdue University,\\
516 Northwestern Avenue, West Lafayette, IN 47907, USA\\
and\\
Weldon School of Biomedical Engineering, Purdue University, \\
206 S. Martin Jischke Drive, West Lafayette, IN 47907, USA\\
and\\
Purdue Center for Cancer Research, Purdue University, \\
201 S. University Street, West Lafayette, IN 47907, USA\\
hectorgomez@purdue.edu}}

\vspace{2mm}

{Guillermo Lorenzo}

{\footnotesize
{Dipartimento di Ingegneria Civile e Architettura, \\
Universit\`{a} degli Studi di Pavia and IMATI-C.N.R.,\\
Via Ferrata 3, 27100 Pavia, Italy \\
guillermo.lorenzo@unipv.it}}

\vspace{2mm}

{Gabriela Marinoschi}

{\footnotesize
{``Gheorghe Mihoc-Caius Iacob'' Institute of Mathematical Statistics \\
and Applied Mathematics of the Romanian Academy,\\
Calea 13 Septembrie 13, 050711 Bucharest, Romania\\
gabriela.marinoschi@acad.ro}}

\vspace{2mm}

{Alessandro Reali}

{\footnotesize
{Dipartimento di Ingegneria Civile e Architettura, \\
Universit\`{a} degli Studi di Pavia and IMATI-C.N.R.,\\
Via Ferrata 3, 27100 Pavia, Italy \\
alessandro.reali@unipv.it}}

\vspace{2mm}

{Elisabetta Rocca}

{\footnotesize
{Dipartimento di Matematica, Universit\`a degli Studi di Pavia and IMATI-C.N.R.,\\
Via Ferrata~5, 27100 Pavia, Italy\\
%and\\
%Istituto di Matematica Applicata e Tecnologie Informatiche ``Enrico Magenes'', CNR,\\
%Via Ferrata~1, 27100 Pavia, Italy\\
elisabetta.rocca@unipv.it
}}

\vspace{5mm}

\end{center}

\begin{center}
{\textbf{Abstract}}
\end{center}

\noindent Chemotherapy is a common treatment for advanced prostate cancer. The standard approach relies on cytotoxic drugs, which aim at inhibiting proliferation and promoting cell death. Advanced prostatic tumors are known to rely on angiogenesis, i.e., the growth of local microvasculature via chemical signaling produced by the tumor. Thus, several clinical studies have been investigating antiangiogenic therapy for advanced prostate cancer, either as monotherapy or in combination with standard cytotoxic protocols. However, the complex genetic alterations that originate and sustain prostate cancer growth complicate the selection of the best chemotherapeutic approach for each patient's tumor. Here, we present a mathematical model of prostate cancer growth and chemotherapy that may enable physicians to test and design personalized chemotherapeutic protocols \emph{in silico}. We use the phase-field method to describe tumor growth, which we assume to be driven by a generic nutrient following reaction-diffusion dynamics. Tumor proliferation and apoptosis (i.e., programmed cell death) can be parameterized with experimentally-determined values. Cytotoxic chemotherapy is included as a term downregulating tumor net proliferation, while antiangiogenic therapy is modeled as a reduction in intratumoral nutrient supply. An additional equation couples the tumor phase field with the production of prostate-specific antigen, which is a prostate cancer biomarker that is extensively used in the clinical management of the disease. We prove the well-posedness of our model and we run a series of representative simulations leveraging an isogeometric method to explore untreated tumor growth as well as the effects of cytotoxic chemotherapy and antiangiogenic therapy, both alone and combined. Our simulations show that our model captures the growth morphologies of prostate cancer as well as common outcomes of cytotoxic and antiangiogenic mono and combined therapy. Additionally, our model also reproduces the usual temporal trends in tumor volume and prostate-specific antigen evolution observed in experimental and clinical studies.    

\vspace{2mm}
{\small
\textbf{Keywords:} prostate cancer; computational oncology; phase field; semilinear parabolic equations; well-posedness; isogeometric analysis

\vspace{2mm}
\textbf{AMS Subject Classification}: 35Q92, 92C50, 65M60, 35K51, 35K58
}

\section{Introduction}

Prostate cancer (PCa) is a major health problem striking over one million men annually and is responsible for over 300,000 deaths \cite{ferlay2015cancer,Mottet2018}. PCa is usually an adenocarcinoma, a form of cancer that originates in the epithelial tissue of the prostate. The evolution of a tumor depends on the genetic alterations that originated it and on its microenvironmental conditions. A key process controlling the conditions of the tumor microenvironment is tumor-induced angiogenesis \cite{figg2008angiogenesis,Hanahan2011}, i.e., the growth of new blood vessels from pre-existing ones via chemical signals produced by the tumor.

Arguably, the most effective way to combat PCa is a combination of prevention and regular screening for early detection. PCa screening is usually accomplished by way of periodic digital rectal exams (DREs) and prostate-specific antigen (PSA) tests \cite{Mottet2018}. The DRE is a physical test in which a doctor palpates the rectal wall next to the prostate to search for hard, lumpy, or abnormal areas typically indicative of cancer. The serum level of PSA is a biomarker of the prostate activity that rises during PCa. The PSA test is a blood test that measures the amount of PSA in the bloodstream. If either the PSA or the DRE test are positive, the patient will be recommended to undergo a biopsy. Together with DRE and serum PSA, the results of the biopsy establish the basis for a diagnosis of PCa. 
Patients diagnosed with localized PCa normally receive a radical treatment with curative intent, such as surgery or radiation \cite{Mottet2018}. 
Advanced PCa patients are usually prescribed androgen-deprivation therapy or chemotherapy \cite{Mottet2018}, which are administered systemically.
In addition to its use in advanced PCa, chemotherapy has also been proposed as an efficient neoadjuvant therapy before radical prostatectomy \cite{Cha2015,pettaway2000neoadjuvant} that aims at reducing the severity of the disease before surgery by shrinking the tumor and eliminating micrometastases that may have developed. 

Chemotherapy for PCa is mostly based on cytotoxic drugs (e.g., docetaxel, cabazitaxel), which obstruct tumor growth by inhibiting cell proliferation and promoting tumor cell death \cite{Mottet2018,Eisenberger2012,Seruga2011}.
While antiangiogenic therapies for PCa are being actively investigated \cite{Antonarakis2012,Seruga2011,Small2012}, they are outside of the standard of care for PCa.
This is somehow counterintuitive because angiogenesis is known to play a central role in the progression of castration-resistant PCa \cite{weidner1993tumor} and an increase in microvascular density is associated to poorer prognosis \cite{mehta2001independent}. There are multiple studies with contradicting evidence about the effectiveness of antiangiogenic therapy \cite{Antonarakis2012,Seruga2011,Small2012,mukherji2013angiogenesis}. However, we would like to highlight that bevacizumab (the most common antiangiogenic drug \cite{Ferrara2004}) did not show significant benefit in castration-resistant PCa when administered alone \cite{reese2001phase}, but it produced a 50\% PSA decline in 75\% of the patients when combined with chemotherapy \cite{picus2011phase}. Thus, bevacizumab and other antiangiogenic drugs are currently being regarded as promising agents to use in combination with cytotoxic chemotherapy. 
However, prostatic tumors have continuously evolving heterogeneous genetic profiles, which may show varying rates of resistance to the prescribed chemotherapeutic treatment and these may even increase during conventional protocols \cite{Seruga2011,Kim2005,Hanahan2011,Gallaher2018}. This is known to be behind the widely varying and sometimes contradictory results of chemotherapeutic clinical studies, also making extremely complex to know whether a patient will benefit from chemotherapy and which drugs are best to treat his tumor.

Recently, the computational modeling and simulation of cancer has shown promise to extend our understanding of these pathologies as well as in forecasting tumor growth and treatment outcomes \cite{Anderson2008,Yankeelov2013,Corwin2013,lorenzo2016tissue,lorenzo2019computer}.
In this context, several studies have been focusing on studying the effects of chemotherapy on tumors through mathematical models and computational simulations \cite{Benzekry2013,Bogdanska2017,Gallaher2018,Hinow2009,Kohandel2007,Powathil2007}. 
This approach would enable to test \emph{in silico} alternative drug protocols and combinations, hence assisting physicians finding optimal chemotherapeutic plans for each patient. Additionally, computational models of cancer growth and treatment may contribute to better comprehend the intricate mechanisms of drug resistance and find early predictors of chemotherapy failure. 

This paper proposes a model of PCa growth and chemotherapy, where the advancement of the tumor is controlled by a critical nutrient. The effect of cytotoxic and antiangiogenic drugs is incorporated by downregulating tumor net proliferation and reducing nutrient supply, respectively. The model also produces as an output the time evolution of serum PSA, which is information commonly used in clinical practice to monitor tumor evolution or the disease's response to a particular treatment. Compared to our previous work on PCa \cite{lorenzo2016tissue,lorenzo2017hierarchically,lorenzo2019computer}, for the model presented herein we perform a complete mathematical analysis, while retaining critical features such as the ability of the model to predict the morphological shift of the tumor under certain circumstances, which is compatible with {\em in vitro} experiments of PCa cell lines \cite{harma2010comprehensive} and clinical observations \cite{Erbersdobler2004,noguchi2000assessment}. We prove that the model is well posed and develop an algorithm to solve the equations numerically. Our computational method is based on isogeometric analysis (IGA), a recently-proposed generalization of the finite element method that uses splines as basis functions \cite{hughes2005isogeometric}. Our computational results show complex tumor dynamics, matching previous observations in both computational and clinical studies.

The paper is organized as follows: Sections 2 and 3.1 describe, respectively, the model equations and the functional framework that we utilize. We prove the well-posedness of the model in Section 3. Section 4 presents our computational method and Section 5 shows representative simulations. We draw conclusions in Section 6.

\section{Model of prostate cancer growth with chemotherapy}
\subsection{Modeling approach}

The model describes the tumor dynamics using a phase field, i.e., a continuous field that defines the time evolution of the tumor's location and geometry. The phase-field method has been extensively used to describe tumor growth in the computational literature \cite{Cavaterra2019,Colli2015,Colli2017,Colli2019,Frigeri2015,Garcke2016,Garcke2018,lorenzo2016tissue,Miranville2019,Frieboes2010,Wise2008,Xu2016}. In this work, the phase field is denoted by $\phi$ and transitions from the value $\phi\approx0$ in the host tissue to $\phi\approx1$ in the tumor. The transition is smooth but steep and takes on a hyperbolic tangent profile in the direction perpendicular to the interface \cite{gomez2018computational}. The model also accounts for the dynamics of a critical nutrient, whose concentration is denoted by $\sigma$ and obeys a reaction-diffusion equation. The concentration of PSA in the prostatic tissue $p$ is governed by a linear reaction-diffusion equation.

\subsection{Governing equations}
\subsubsection{Tumor dynamics}\label{tumoreq}
\begin{figure}[t]
\centerline{\includegraphics[width=\linewidth]{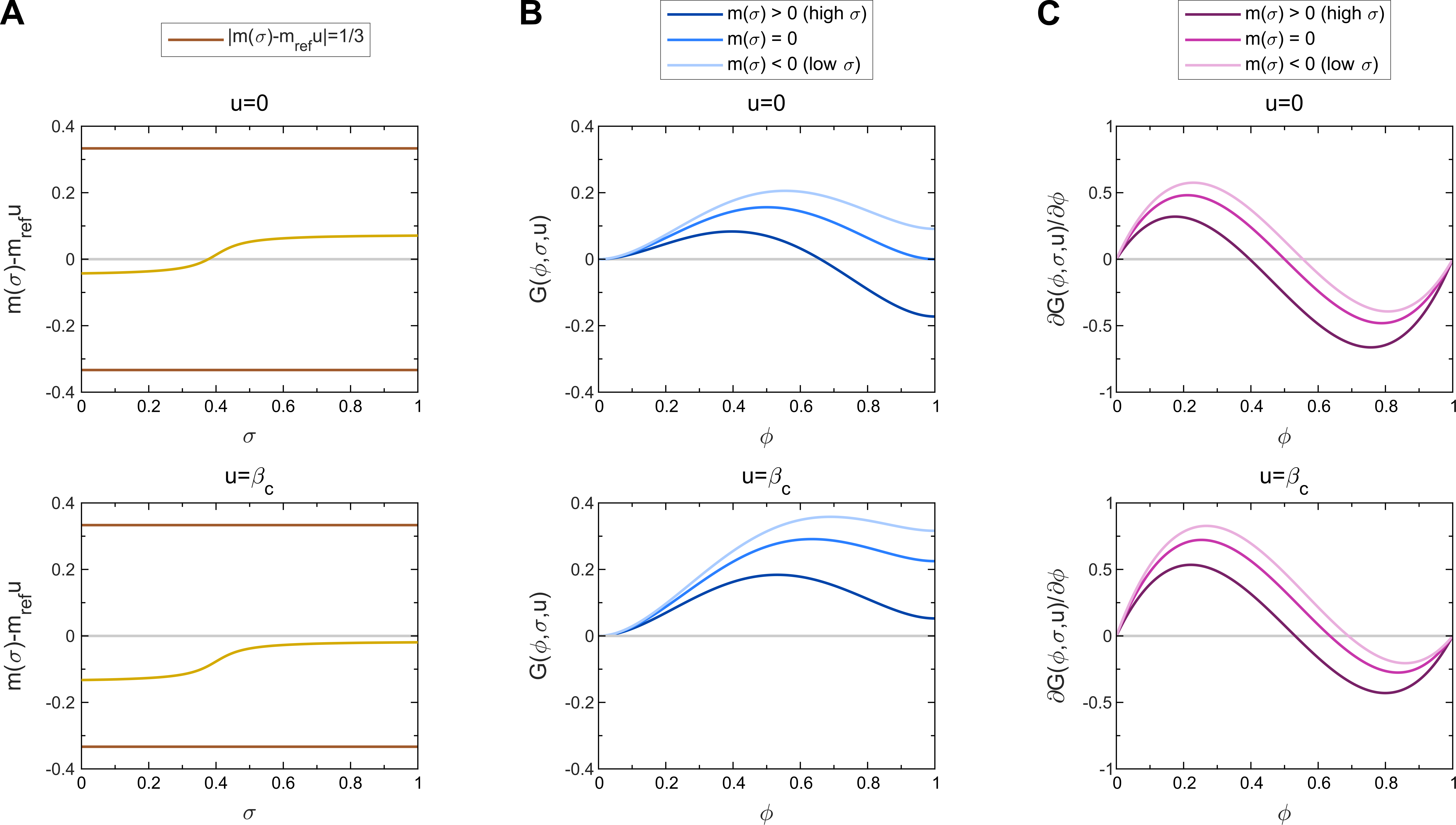}}
\vspace*{8pt}
\caption{A high nutrient environment energetically favors tumor growth in our model, whereas low nutrient availability and the action of the cytotoxic drug energetically obstruct it. 
(A) Plot of $m(\sigma)-m_{ref}u$ with respect to the values of $\sigma$ for $u=0$ (top) and $u=\beta_c$ (bottom). 
(B) Plot of $G\left(\phi,\sigma,u\right)$ with respect to the values of $\phi$ for three values of function $m(\sigma)$ (positive, zero, and negative) combined with $u=0$ (top) and $u=\beta_c$ (bottom).
(C) Plot of $\partial G\left(\phi,\sigma,u\right)/\partial\phi$ with respect to the values of $\phi$ for three values of function $m(\sigma)$ (positive, zero, and negative) combined with $u=0$ (top) and $u=\beta_c$ (bottom).}
\label{functionplots}
\end{figure}

The tumor dynamics is described by the equation \cite{Xu2016}
\begin{equation}
\phi _{t}=\lambda \Delta \phi -\frac{\partial G}{\partial \phi }(\phi
,\sigma ,u) 
\end{equation}
where $\phi_t$ denotes the partial derivative of $\phi$ with respect to time, $\Delta$ is the Laplace operator, $\lambda$ is the diffusion coefficient of tumor cells and
\begin{equation}\label{gfunct}
G(\phi,\sigma,u)=F(\phi)-h(\phi)(m(\sigma)-m_{ref}u).
\end{equation}
The diffusion coefficient of tumor cells can be computed as $\lambda=M\ell^2$, where both $M$ and $\ell$ are positive real constants denoting the tumor mobility and interface width, respectively \cite{Xu2016}.
Here, $F(\phi)=M\phi^2(1-\phi)^2$ is a double-well potential, i.e., a nonconvex function, typical in phase-field modeling, which allows the coexistence of the tumoral ($\phi\approx1$) and healthy ($\phi\approx0$) tissue. 
%The constant $M$ takes on positive values. 
The function $h(\phi)=M\phi^2(3-2\phi)$ is also common in non-conserved dynamics of phase fields. It is usually called interpolation function because it verifies the properties $h(0)=0$ and $h(1)=1$; another important property of $h$ is that $h^\prime(0)=h^\prime(1)=0$, where $h^\prime$ denotes the derivative of $h$. The function $m(\sigma)$ is normally called tilting function and, in our model, is defined as
\begin{equation}\label{msigma}
m(\sigma )=m_{ref}\left( \frac{\rho +A}{2}+\frac{\rho -A}{\pi }\arctan \left( \frac{\sigma -\sigma _{l}}{\sigma _{r}}\right) \right),
\end{equation}
where $\rho$ and $A$ represent constant proliferation and apoptosis indices, respectively. These nondimensional parameters are related  to the proliferation and apoptotis rates in tumoral tissue as follows:
\begin{equation}\label{rhodef}
\rho=\frac{K_\rho}{\bar{K_\rho}} 
\end{equation}
and
\begin{equation}\label{Adef}
A=-\frac{K_A}{\bar{K_A}},
\end{equation}
 where $K_\rho$ is the proliferation rate of tumor cells, $K_A$ is the apoptosis rate of tumor cells, $\bar{K_\rho}$ is a scaling reference value for the proliferation rate, and $\bar{K_A}$  is a scaling reference value for the apoptosis rate. Therefore, $m(\sigma)$ can be interpreted as a function describing the tumor net proliferation rate. The positive constant $m_{ref}$ scales the strength of the tilting function within our phase-field framework. The constants $\sigma_r$ and $\sigma_l$ are, respectively, a reference and a threshold value for the nutrient concentration \cite{Xu2016}. For nutrient concentrations lower than $\sigma_l$ healthy tissue is energetically more favorable than tumor tissue and viceversa. The function $u$ in \eqref{gfunct} represents the tumor-inhibiting effect of a cytotoxic drug and is described in Section \ref{chemotherapy}. When $|m(\sigma)-m_{ref}u|<1/3$, the function $G$ is a double-well potential with local minima at $\phi=0$ and $\phi=1$. Within this range, low values of the nutrient concentration (or large values of $u$) produce a lower energy level (value of $G$) in the healthy tissue ($\phi=0$) than in the tumoral tissue ($\phi=1$). The opposite is true for high values of the nutrient concentration (or low values of $u$). The described behavior of  $G$ is further illustrated in Figure \ref{functionplots}. From \eqref{gfunct}, we can obtain
\begin{equation}
\frac{\partial G}{\partial\phi}=2\phi (1-\phi )f(\phi,\sigma ,u) 
\end{equation}
where%
\begin{equation}\label{f}
f(\phi ,\sigma ,u)=M\left[ 1-2\phi -3\left( m(\sigma)-m_{ref}u\right) \right].  
\end{equation}

\subsubsection{Cytotoxic chemotherapy} \label{chemotherapy}

Docetaxel is a usual cytotoxic drug in the clinical management of advanced PCa \cite{Mottet2018,Seruga2011,Eisenberger2012}. The standard cytotoxic chemotherapy based on docetaxel usually consists of up to 10 drug doses delivered every three weeks.  We consider that the action of the cytotoxic drug on tumor dynamics depends linearly on the drug concentration
 \cite{Hahnfeldt2003,Colli2019,Garcke2018,Bogdanska2017,Bodnar2019,Kohandel2007,Hinow2009,Powathil2007}. The pharmacodynamics of this cytotoxic drug shows an approximately exponential decrease in drug concentration after delivery of the dose at the systemic level \cite{Baker2004,Tije2005}. Hence, we propose the formulation
\begin{equation}\label{eq_u}
u(t)=\sum_{i=1}^{n_c}\beta_c d_c e^{-\frac{t-T_{c,i}}{\tau_c}}\mathcal{H}\left(t-T_{c,i}\right),
\end{equation}
where $n_c$ is the number of chemotherapy cycles, $\beta_c$ measures the cytotoxic effect of the treatment on tumor dynamics per unit of drug dose delivered, $d_c$ is the drug dose,  $T_{c,i}$ with $i=1,\dots,n_c$ are the times of drug delivery, $\tau_c$ is the mean lifetime of the chemotherapeutic drug, and $\mathcal{H}$ denotes the Heaviside function. These parameters are required to be strictly positive. The experimental and clinical values observed for $\beta_c$ ensure that the condition $|m(\sigma)-m_{ref}u|<1/3$ is verified; see Section \ref{tumoreq}. According to \eqref{gfunct}, $u(t)$ can be interpreted as a cytotoxic drug-induced decrease in tumor net proliferation that decays as the drug concentration in time. Thus, the action of the cytotoxic drug energetically disfavors tumor tissue (see Figure \ref{functionplots}). Notice that the formulation paradigm in \eqref{eq_u} is very similar to the modeling of radiation effects on tumor dynamics \cite{Lima2017,Perez-Garcia2015,Corwin2013}. In principle, the formulation proposed in \eqref{eq_u} is extensible to other cytotoxic drugs by adapting $\tau_c$, $\beta_c$ and $d_c$ accordingly.

\subsubsection{Nutrient dynamics}

Nutrient dynamics is controlled by the equation
\begin{equation}\label{eq_nut}
\sigma_t=\eta\Delta\sigma+S_h(1-\phi)+(S_c-s)\phi-\left( \gamma_h(1-\phi) + \gamma_c\phi \right) \sigma,
\end{equation}
where $\eta$ is the diffusion coefficient of the nutrient, $S_h$ is the nutrient supply rate in the healthy tissue, $S_c$ stands for the nutrient supply rate in the cancerous tissue, $s$ is a given function yielding the reduction in nutrient supply caused by antiangiogenic therapy (see Section~\ref{antiangiogenic}), and $\gamma_h$, $\gamma_c$ are positive constants that represent the nutrient uptake rate in the healthy and cancerous tissue, respectively.  We require that $S_h$, $S_c$, and $s$ are non-negative and $s$ should satisfy the constraint $s\leq S_c$.

%The term 
%\begin{equation}\label{nsupply}
%S=S_h(1-\phi)+S_c\phi-v\phi
%\end{equation}
%represents the nutrient supply rate in the prostatic tissue. The coefficient $S_h$ is the nutrient supply rate in the healthy tissue, while $S_c$ stands for the nutrient supply rate in the cancerous tissue. The function $v$ represents the reduction in nutrient supply caused by antiangiogenic therapy, whose formulation we describe below. We require that $S_h$, $S_c$, and $v$ are non-negative and $v$ should satisfy the constraint $v\leq S_c$.

\subsubsection{Antiangiogenic therapy}\label{antiangiogenic}

Bevacizumab is a common antiangiogenic drug that has been actively investigated for the treatment of PCa \cite{Antonarakis2012,Seruga2011,Small2012,Ferrara2004}. This antiangiogenic therapy is usually delivered simultaneously with the cytotoxic chemotherapy or following a similar schedule in case it is the only treatment.
 Again, we consider that the action of this antiangiogenic drug depends linearly on its concentration \cite{Kohandel2007,Hahnfeldt1999,Benzekry2013,Benzekry2012}. Pharmacodynamic studies of bevacizumab also revealed an exponential decay in drug concentration following systemic delivery of the prescribed dose \cite{Gordon2001,Lu2008}. Therefore, we choose an analogous formulation to \eqref{eq_u}:
\begin{equation}\label{eq_v}
s(t)=\sum_{i=1}^{n_a}\beta_a d_a e^{-\frac{t-T_{a,i}}{\tau_a}}\mathcal{H}\left(t-T_{a,i}\right)
\end{equation}
where $n_a$ is the number of antiangiogenic treatment cycles, $\beta_a$ measures the antiangiogenic effect of the treatment on the nutrient supply per unit of drug dose delivered, $d_a$ is the drug dose,  $T_{a,i}$ are the times of drug delivery ($i=1,\dots,n_a$), $\tau_a$ stands for the mean lifetime of the antiangiogenic drug, and $\mathcal{H}$ denotes the Heaviside function. These parameters are required to be strictly positive. Again, the formulation in \eqref{eq_v} may be extended to other antiangiogenic drugs by appropriately calibrating $\tau_a$, $\beta_a$ and $d_a$.

\subsubsection{Tissue PSA dynamics}

Both healthy and cancerous prostatic cells release PSA, although tumor cells do so at a much larger rate. The PSA is assumed to diffuse through the prostatic tissue and decay naturally at rate $\gamma_p$. Hence, we propose the equation \cite{lorenzo2016tissue} 
\begin{equation}\label{tpsa}
p_t = D\Delta p + \alpha_h(1-\phi) + \alpha_c \phi - \gamma_p p
\end{equation}
where $D$ is the diffusion constant. The constants $\alpha_h$ and $\alpha_c$ represent, respectively, the tissue PSA production rate of healthy and malignant cells. The serum PSA $P_s$ can be defined as the integral of the tissue PSA over the prostate $\Omega$, that is,
\begin{equation}
P_s=\int_\Omega p \, dx.
\end{equation}
By using \eqref{tpsa}, one can show that for free-flux boundary conditions
\begin{equation}\label{spsa}
\frac{d P_s}{d t}=\alpha_h V_h + \alpha_c V_c - \gamma_p P_s,
\end{equation}
where $V_h$ is the volume of prostatic healthy tissue and $V_c$ is the tumor volume. Equation \eqref{spsa}, which we derived from the tissue level equation \eqref{tpsa}, was proposed as a phenomenological model of serum PSA dynamics that fits clinical data \cite{swanson2001quantitative}.

\section{Well-posedness of the model}\label{wellpos}
\setcounter{equation}{0}

In this section we shall provide an existence and uniqueness result for this
system 
\begin{equation}
\phi _{t}-\lambda \Delta \phi +2\phi (1-\phi )f(\phi ,\sigma ,u)=0,\mbox{ in 
}\piergab{Q_T :={}} (0,T)\times \Omega ,  \label{e1}
\end{equation}%
\begin{equation}
\sigma _{t}-\eta \Delta \sigma +\gamma _{h}\sigma +(\gamma _{c}-\gamma
_{h})\sigma \phi \pcol{{}= \piergab{S_h + (S_c -S_h)\phi -s \phi ,{}}\mbox{ in }Q_T} ,  \label{e2}
\end{equation}%
\begin{equation}
p_{t}-D\Delta p+\gamma _{p}p=\alpha _{h}+(\alpha _{c}-\alpha _{h})\phi ,%
\mbox{ in }\piergab{Q_T} ,  \label{e3}
\end{equation}%
\begin{equation}
\phi =0,\mbox{ }\frac{\partial \sigma }{\partial \nu }=\frac{\partial p}{%
\partial \nu }=0,\mbox{ on }\piergab{ \Sigma_T:={}} (0,T)\times \partial \Omega ,  \label{e4}
\end{equation}%
\begin{equation}
\phi (0)=\phi _{0},\mbox{ }\sigma (0)=\sigma _{0},\mbox{ }p(0)=p_{0},\mbox{
in }\Omega ,  \label{e5}
\end{equation}%
where $f$ is given by (\ref{f}) and (\ref{msigma}); $M, m_{ref}, \rho, A, \sigma_l, \sigma_r$ are given real constants; and  $\lambda ,$ $\eta ,$ $D,$ $\gamma
_{c},$ $\gamma _{h},$ $\gamma _{p},$ $\alpha _{c},$ $\alpha _{h}$, \piergab{$S_c$, $S_h$} are positive fixed real constants.

Here, $\Omega $ is an open bounded subset of \pier{$\mathbb{R}^{N}$, $N\leq3$,} with a
sufficiently smooth boundary $\partial \Omega ,$ $\frac{\partial }{\partial
\nu }$ represents the outward normal derivative to $\partial \Omega $ and \BBB $T$
denotes a final time. \EEE

\subsection{Functional framework}

The problem will be treated in a functional framework involving the space $%
H=L^{2}(\Omega ),$ which is identified with its dual space $H^{\prime
}\pier{\cong} H,$ and the Sobolev spaces%
\[
V_{0}=H_{0}^{1}(\Omega ),\mbox{ }V_{0}^{\prime }=(H_{0}^{1}(\Omega
))^{\prime }=:H^{-1}(\Omega ),\mbox{ \ }V=H^{1}(\Omega ),\mbox{ }V^{\prime
}=(H^{1}(\Omega ))^{\prime },\mbox{ } 
\]%
\[
W_{0}=H^{2}(\Omega )\cap H_{0}^{1}(\Omega ),\mbox{ \ }W=\left\{ y\in
H^{2}(\Omega );\mbox{ }\frac{\partial y}{\partial \nu }=0\mbox{ on }\partial
\Omega \right\} , 
\]%
and with the dense and compact injections $W_{0}\subset V_{0}\subset H\subset
V_{0}^{\prime }$ and $W\subset V\subset H\subset V^{\prime }.$ 
$H^1_0(\Omega)$ contains the elements of $H^1(\Omega)$ with null trace on the boundary $\partial\Omega$.
\piergab{For the norms in these spaces we will use the notation $\| \, \cdot \, \|_B$, where 
$B$ is the space we are considering.}

\piergab{For \BBB $q\in \lbrack 1,\infty ]$ and $z\in L^{q}(\Omega )$ or $z\in
L^{q}(Q_{T}) $, we simply denote \pier{the norm of $z$ in these $L^q$ spaces} by $\left\Vert z\right\Vert
_{q}.$\EEE}

\piergab{We will also 
make use of spaces of functions that depend on time with values in a Banach space $B$. Namely, for 
 $q\in \lbrack 1,\infty ]$ 
we consider the space $L^q (0,T;B)$ of measurable functions $t\mapsto z(t)$ such that 
$t \mapsto \| z(t)\|_B^q $ is integrable on $(0,T)$ (or essentially bounded  if $q=\infty$) 
and the space $C([0,T]; B) $ of continuous functions from $[0,T]$ to $B$. Perhaps, it is important to 
note  that $L^2 (0,T;H)$ is a space completely isomorph to $L^2 (Q_T)$. Moreover, for 
 $q\in \lbrack 1,\infty ]$,  $W^{1,q} (0,T;B)$
will denote the space of functions $t\mapsto z(t)$ such that both $z$ and its (weak) derivative $z_t$
belong to $L^q (0,T;B)$. We point out that  $W^{1,q} (0,T;B) \subset  C([0,T]; B) $
for all  $q\in \lbrack 1,\infty ]$.} 
\medskip

\textbf{Hypotheses. }Following the considerations presented in the
Introduction \piergab{and extending the particular choices of $u$ and 
$s$ in \eqref{eq_u} and \eqref{eq_v}, we let}
\begin{equation}
\pier{u\in L^{\infty }(Q_{T}),\mbox{ }\ \piergab{s\in L^{\infty }(Q_{T})},}  \label{e6}
\end{equation}%
and it can be easily checked that 
\begin{equation}
(r_{1},r_{2},r_{3})\rightarrow f(r_{1},r_{2},r_{3})\mbox{ is globally
Lipschitz continuous on }\mathbb{R}^{3}.  \label{e7-0}
\end{equation}%
We also set 
\begin{equation}
\gamma _{ch}=\gamma _{c}-\gamma _{h},\quad \alpha _{ch}=\alpha _{c}-\alpha_{h}, 
\quad \piergab{S_{ch}= S_c - S_h},  \label{e7-2}
\end{equation}%
and point out that the problem parameters $\lambda ,$ $\eta ,$ $D,$ $\gamma
_{c},$ $\gamma _{h},$ $\gamma_{p},$ $\alpha _{c},$ $\alpha _{h}$, \piergab{$S_c$, $S_h$} are positive.

\medskip

\piergab{In the sequel, by $C$ we shall denote a constant, that may
change from line to line, depending on the problem parameters, the domain $%
\Omega ,$ the final time $T,$\ the norms of the initial data and possibly on
the norms of $u$ and $s$. Moreover, we assume that (\ref{e6})-(\ref{e7-2}) hold.}

%We set $Q_{T}:=(0,T)\times \Omega ,$ $\Sigma _{T}:=(0,T)\times \partial \Omega .$

\medskip

\noindent \textbf{Definition 3.1.} Let $(\varphi _{0},\sigma _{0},p_{0})\in
H\times H\times H.$ \pier{A solution to the system (\ref{e1})-(\ref{e5}) is a triplet $(\phi,\sigma,p)$, with}
\begin{eqnarray}
\phi &\in &W^{1,2}(0,T;V_{0}^{\prime })\cap C([0,T];H)\cap
L^{2}(0,T;V_{0})\cap L^{\infty }(Q_{T}), \label{e8} \\
\sigma &\in &W^{1,2}(0,T;V^{\prime })\cap C([0,T];H)\cap L^{2}(0,T;V),
 \nonumber \\
p &\in &W^{1,2}(0,T;V^{\prime })\cap C([0,T];H)\cap L^{2}(0,T;V),  \nonumber
\end{eqnarray}%
which satisfies 
\begin{align}
&\int_{0}^{T}\left\langle \phi _{t}(t),\psi _{1}(t)\right\rangle
_{\pcol{V_0^{\prime },V_0}}dt+\int_{Q_{T}}\left\{ \lambda \nabla \phi \cdot \nabla \psi
_{1}+2\phi (1-\phi )f(\phi ,\sigma ,u)\psi _{1}\right\} dxdt
\label{e9} \\
&+\int_{0}^{T}\left\langle \sigma _{t}(t),\psi _{2}(t)\right\rangle
_{V^{\prime },V}dt+\int_{Q_{T}}(\eta \nabla \sigma \cdot \nabla \psi
_{2}+\gamma _{h}\sigma \psi _{2}+\gamma _{ch}\sigma \phi \psi _{2})dxdt 
\nonumber \\
&=\int_{Q_{T}}\piergab{(S_h\psi _{2} + (S_{ch} - s) \phi \psi_2)}dxdt,\mbox{ for all }(\psi _{1},\psi _{2})\in
L^{2}(0,T;V_{0}\times V),  \nonumber
\end{align}%
\begin{eqnarray}
&&\int_{0}^{T}\left\langle p_{t}(t),\psi (t)\right\rangle _{V^{\prime
},V}dt+\int_{Q_{T}}(D\nabla p \cdot \nabla \psi +\gamma _{p}p\psi )dxdt
\label{e10} \\
&&=\int_{Q_{T}}(\alpha _{h}+\pier{\alpha _{ch}}\phi )\psi dxdt,\mbox{
for all }\psi \in L^{2}(0,T;V),  \nonumber
\end{eqnarray}%
and 
\begin{equation}
(\varphi ,\sigma ,p)(0)=(\varphi _{0},\sigma _{0},p_{0}).  \label{e10-0}
\end{equation}
\pier{Please note that taking  in (\ref{e9}) $\psi_2 $ first, and $\psi_1 $ then, equal to the null function
(and this is of course a suitable choice in both cases), we obtain two separate variational equalities for 
$\phi$ and $\sigma$, as it is (\ref{e10}) for $p$.}
\medskip

Let us denote%
\begin{eqnarray}
X_{0} &=&W^{1,2}(0,T;H)\cap \pier{C([0,T];}V_{0})\cap L^{2}(0,T;W_{0}),
\label{e11} \\
X &=&W^{1,2}(0,T;H)\cap \pier{C([0,T];}V)\cap L^{2}(0,T;W).  \nonumber
\end{eqnarray}

\medskip

\noindent \textbf{Theorem 3.2. }\textit{Let } \textit{\ }%
\begin{equation}
(\phi _{0},\sigma _{0},p_{0})\in H\times H\times H,  \label{e12}
\end{equation}%
\begin{equation}
0\leq \phi _{0}(x)\leq 1\mbox{ \textit{a.e.} }x\in \Omega .  \label{e12-0}
\end{equation}%
\textit{Then, \pier{the} system} (\ref{e1})-(\ref{e5}) \textit{has a unique solution \pier{$(\phi,\sigma,p)$} in
the sense of }Definition 3.1, \textit{such that} 
\[
0\leq \phi (t,x)\leq 1\mbox{ \textit{a.e.} }(t,x)\in Q_{T}. 
\]%
\textit{If }$(\sigma _{0},p_{0})\in L^{\infty }(\Omega )\times L^{\infty
}(\Omega ),$\textit{\ then} $(\sigma ,p)\in L^{\infty }(Q_{T})\times
L^{\infty }(Q_{T}).$ \textit{Moreover, if} 
\begin{equation}
\sigma _{0}(x)\geq 0,\mbox{ }p_{0}(x)\geq 0\mbox{ \ \textit{a.e.} }x\in \Omega
,\mbox{ \ \ }\piergab{s(t,x)\leq S_c}\mbox{ \ \ \textit{a.e.} }(t,x)\in Q_{T}.  \label{e13-0}
\end{equation}%
\textit{then we have} 
\begin{equation}
\sigma (t,x)\geq 0,\mbox{ }p(t,x)\geq 0\mbox{ \textit{a.e.} }(t,x)\in Q_{T}.
\label{e13-2}
\end{equation}%
\textit{Finally, if }$(\phi _{0},\sigma _{0},p_{0})\in V_{0}\times V\times
V, $\textit{\ the solution has the supplementary regularity }$(\phi ,\sigma
,p)\in X_{0}\times X\times X,$ \textit{and satisfies the estimate }%
\begin{eqnarray}
&&\left\Vert \phi \right\Vert _{X_{0}}+\left\Vert \sigma \right\Vert
_{X}+\left\Vert p\right\Vert _{X}  \label{e13} \\
&\leq &C(\left\Vert \phi _{0}\right\Vert _{V_{0}}^{2}+\left\Vert \sigma
_{0}\right\Vert _{V}^{2}+\left\Vert p_{0}\right\Vert _{V}^{2}+\left\Vert
u\right\Vert _{\pier{L^{2}(0,T;H)}}^{2}+\left\Vert \piergab{s}\right\Vert _{\pier{L^{2}
(0,T;H)}}^{2}+1).  \nonumber
\end{eqnarray}%
\textit{In addition, the solution is continuous with respect to the data,
that is, for two solutions }$(\phi _{i},\sigma _{i},p_{i})$\textit{\
corresponding to \pier{the data}} $(\phi _{0}^{i},\sigma _{0}^{i},p_{0}^{i},u_{i}, \piergab{s}_{i}),$%
\textit{\ }$i=1,2,$\textit{\ we have} 
\begin{eqnarray}
&&\left\Vert (\phi _{1}-\phi _{2})(t)\right\Vert _{H}^{2}+\left\Vert (\sigma
_{1}-\sigma _{2})(t)\right\Vert _{H}^{2}+\left\Vert
(p_{1}-p_{2})(t)\right\Vert _{H}^{2}  \label{e14} \\
&&+\left\Vert \phi _{1}-\phi _{2}\right\Vert
_{L^{2}(0,T;V_{0})}^{2}+\left\Vert \sigma _{1}-\sigma _{2}\right\Vert
_{L^{2}(0,T;V)}^{2}+\left\Vert p_{1}-p_{2}\right\Vert _{L^{2}(0,T;V)}^{2} 
\nonumber \\
&\leq &C\left( \left\Vert \phi _{0}^{1}-\phi _{0}^{2}\right\Vert
_{H}^{2}+\left\Vert \sigma _{0}^{1}-\sigma _{0}^{2}\right\Vert
_{H}^{2}+\left\Vert p_{0}^{1}-p_{0}^{2}\right\Vert _{H}^{2}\right.  \nonumber
\\
&&\left. \mbox{ \ \ \ \ }\pcol{+\left\Vert u_{1}-u_{2}\right\Vert
_{L^{2}(0,T;H)}^{2}+\left\Vert \piergab{s}_{1}- \piergab{s}_{2}\right\Vert
_{L^{2}(0,T;H)}^{2}}\right)  \nonumber
\end{eqnarray}%
\textit{for all} $t\in \lbrack 0,T].$

\textit{\medskip }

\noindent \textbf{Proof. }We shall prove the existence of the solution by
using the Banach fixed point theorem combined with a variational approach.
\BBB In this respect, we will \EEE be formal in the following sense: when referring 
to the weak solutions of equations like  (\ref{e2})-(\ref{e4}) coupled with the related 
boundary conditions, let us write directly the equations instead of their variational formulations. 
This would allow us to come quickly to the point in our argumentations.

Let us set 
\[
\mathcal{M}=\{z\in C([0,T];H)\pier{:}\mbox{ }0\leq z(t,x)\leq 1\mbox{ a.e. }%
(t,x)\in Q_{T}\}, 
\]%
which is a complete metric space \pier{provided we take the distance induced by some norm in $C([0,T];H)$.}

Let us fix $z\in \mathcal{M}$ in equations (\ref{e1}), (\ref{e2}):%
\begin{equation}
\phi _{t}-\lambda \Delta \phi =-2\phi (1-\phi )f(z,\sigma ^{z},u),\mbox{ in }%
Q_{T},  \label{e15-0}
\end{equation}%
\begin{equation}
\sigma _{t}-\eta \Delta \sigma +\gamma _{h}\sigma +\gamma _{ch}\sigma z= \piergab{S_h+ (S_{ch} - s ) z},%
\mbox{ in }Q_{T},  \label{e15}
\end{equation}%
where $\sigma ^{z}$ in (\ref{e15-0}) is the solution to (\ref{e15}),
corresponding to $z$, with $\frac{\partial \sigma }{\partial \nu }=0$ on $%
\Sigma _{T}$ and $\sigma (0)=\sigma _{0}$ in $\Omega .$ Thus, we can define
the mapping 
\[
\pier{z\mapsto  \phi ^{z}}=:\Psi (z) 
\]%
where $\phi ^{z}$ is the solution to (\ref{e15-0}), corresponding to $z$ \pier{and $\sigma^z$},
with\ the homogeneous Dirichlet boundary condition and the initial datum $%
\phi _{0}.$ We shall prove that $\Psi (\mathcal{M})\subset \mathcal{M}$ and
that $\Psi $ is a contraction mapping on $\mathcal{M}$.

First, we \pier{treat the initial-boundary value problem for equation \piergab{(\ref{e15})}.} For all $t\in \lbrack 0,T]$ we
introduce the operator $A(t):V\rightarrow V^{\prime }$ \pier{defined} by 
\[
\left\langle A(t)\sigma ,\psi \right\rangle _{V^{\prime },V}=\int_{\Omega
}(\eta \nabla \sigma \cdot \nabla \psi +\gamma _{h}\sigma \psi +\gamma
_{ch}z\sigma \psi )dx,\pier{\mbox{ for } \sigma, \psi} \in V, 
\]%
and observe that it is continuous 
\[
\left\Vert A(t)\sigma \right\Vert _{V^{\prime }}\leq \max \{\eta ,\gamma
_{h},\left\vert \gamma _{ch}\right\vert \}\left\Vert \sigma \right\Vert _{V} 
\]%
and quasi-monotone from $V$ to $V^{\prime }$, i.e.,%
\[
\left\langle A(t)\sigma ,\sigma \right\rangle _{V^{\prime },V}\geq \min
\{\eta ,\gamma _{h}\}\left\Vert \sigma \right\Vert _{V}^{2}-\left\vert
\gamma _{ch}\right\vert \left\Vert \sigma \right\Vert _{H}^{2}. 
\]%
Due to \piergab{a general solvability result  for linear parabolic problems (the reader 
may consult, e.g., Ref.~ \cite{Lions1961}}), 
for $\sigma _{0}\in H$ and $\piergab{s} \in L^{\infty }(Q_{T})$ the
initial-boundary value problem for (\ref{e15}) has a unique solution 
\[
\sigma ^{z}\in W^{1,2}(0,T;V^{\prime })\cap C([0,T];H)\cap
L^{2}(0,T;V). 
\]%
A first estimate performed by testing (\ref{e15}) by $\sigma ^{z}$ and
integrating over $(0,t)$ yields, after a standard calculation,%
\begin{equation}
\left\Vert \sigma ^{z}(t)\right\Vert _{H}^{2}+\int_{0}^{t}\left\Vert \sigma
^{z}(\tau )\right\Vert _{V}^{2}d\tau \leq C\left(\piergab{1 +{}}\left\Vert \sigma
_{0}\right\Vert _{H}^{2}+\left\Vert \piergab{s} \right\Vert _{L^{2}(0,T;H)}^{2}\right),  \label{e16}
\end{equation}%
for all $t\in \lbrack 0,T],$ where we used the fact that $z\in \mathcal{M}.$ As specified before, $C$ is
a constant, that may vary from line to line, depending on the problem
parameters and $T$\pier{: indeed, the Gronwall lemma has been applied to derive 
(\ref{e16}).}

Now, if we \pier{take} two solutions $\sigma _{1}$ and $\sigma _{2}$
corresponding to \pcol{$(z_{1},\piergab{s_{1}} ,\sigma_0^1)$ and $(z_{2},\piergab{s_{2}},\sigma_0^2),$} 
\pier{thanks to the H\"older and Young inequalities}
we obtain for the differences $\sigma :=\sigma _{1}-\sigma _{2}$
and $z:=z_{1}-z_{2}:$ 
\begin{align}
&\frac{1}{2}\left\Vert \sigma (t)\right\Vert _{H}^{2}+\pier{\min} \{\eta ,\gamma
_{h}\}\int_{0}^{t}\left\Vert \sigma (\tau )\right\Vert _{V}^{2}d\tau
\label{e16-0} \\
&\leq \frac{1}{2}\left\Vert \pcol{\sigma _{0}^1-\sigma_0^2}\right\Vert
_{H}^{2}+\left\vert \gamma _{ch}\right\vert \int_{0}^{t}\int_{\Omega
}\left\{ \left\vert z\sigma _{1}\right\vert +\left\vert \sigma
z_{2}\right\vert \right\}\pier{\vert\sigma\vert} dxd\tau \nonumber \\
&\quad{} 
\piergab{{}+\int_{0}^{t}\int_{\Omega }\left\vert S_{ch}z\sigma\right\vert dxd\tau 
+\int_{0}^{t}\int_{\Omega }\left(\vert s_1 z\vert + | (s_1 - s_2) z_2| \right)|\sigma\vert  dxd\tau {}}
    \nonumber \\ 
&\leq \frac{1}{2}\left\Vert\pcol{\sigma _{0}^1-\sigma_0^2}\right\Vert
_{H}^{2}+\left\vert \gamma _{ch}\right\vert \int_{0}^{t}\left\Vert z(\tau
)\right\Vert _{H}\left\Vert \sigma _{1}(\tau )\right\Vert _{4}\left\Vert
\sigma (\tau )\right\Vert _{4}d\tau   \nonumber \\
&\quad{}+\left\vert \gamma _{ch}\right\vert
\int_{0}^{t}\left\Vert \sigma (\tau )\right\Vert _{H}^{2}d\tau
\piergab{{}+ |S_{ch}|^2 \int_{0}^{t}\left\Vert z (\tau )\right\Vert _{H}^{2}d\tau
+ 2 \int_{0}^{t} \| s_1\|_\infty^2 \left\Vert z (\tau )\right\Vert _{H}^{2}d\tau}
\nonumber \\
&\quad{}
\piergab{{}+2\int_{0}^{t}\left\Vert (s_{1}-s_{2})(\tau )\right\Vert _{H}^{2}d\tau
+\frac{1}{2}\int_{0}^{t}\left\Vert \sigma (\tau )\right\Vert _{H}^{2}d\tau . } \nonumber
\end{align}%
By the Sobolev embedding theorems, it turns out that $V\subset L^{4}(\Omega )$
with continuous embedding. Then, by exploiting the Young inequality and (\ref%
{e16}) \pier{for $\sigma_1$} we infer that 
\begin{eqnarray*}
&&\left\vert \gamma _{ch}\right\vert \int_{0}^{t}\left\Vert z(\tau
)\right\Vert _{H}\left\Vert \sigma _{1}(\tau )\right\Vert _{4}\left\Vert
\sigma (\tau )\right\Vert _{4}d\tau \leq C\int_{0}^{t}\left\Vert z(\tau
)\right\Vert _{H}\left\Vert \sigma _{1}(\tau )\right\Vert _{V}\left\Vert
\sigma (\tau )\right\Vert _{V}d\tau \\
&\leq &\frac{1}{2}\pier{\min} \{\eta ,\gamma _{h}\}\int_{0}^{t}\left\Vert \sigma
(\tau )\right\Vert _{V}^{2}d\tau +C\int_{0}^{t}\left\Vert z(\tau
)\right\Vert _{H}^{2}\left\Vert \sigma _{1}(\tau )\right\Vert _{V}^{2}d\tau
\\
&\leq &\frac{1}{2}\pier{\min} \{\eta ,\gamma _{h}\}\int_{0}^{t}\left\Vert \sigma
(\tau )\right\Vert _{V}^{2}d\tau +C\left\Vert z\right\Vert
_{C([0,t];H)}^{2}\int_{0}^{t}\left\Vert \sigma _{1}(\tau )\right\Vert
_{V}^{2}d\tau \\
&\leq &\frac{1}{2}\pier{\min} \{\eta ,\gamma _{h}\}\int_{0}^{t}\left\Vert \sigma
(\tau )\right\Vert _{V}^{2}d\tau +C\left\Vert z\right\Vert
_{C([0,t];H)}^{2}\left( \piergab{1+\left\Vert \pcol{\sigma _{0}^1}\right\Vert _{H}^{2}+\left\Vert
s_{1}\right\Vert _{L^{2}(0,T;H)}^{2}}\right).
\end{eqnarray*}%
Going back to (\ref{e16-0}) we deduce that%
\begin{eqnarray*}
&&\left\Vert \sigma (t)\right\Vert _{H}^{2}+\pier{\min} \{\eta ,\gamma
_{h}\}\int_{0}^{t}\left\Vert \sigma (\tau )\right\Vert _{V}^{2}d\tau \\
&\leq &\left\Vert\pcol{\sigma _{0}^1-\sigma_0^2}\right\Vert
_{H}^{2}+4\int_{0}^{t}\left\Vert (\piergab{s_{1}-s_{2}})(\tau )\right\Vert
_{H}^{2}d\tau +C_{01}\left\Vert z\right\Vert
_{C([0,t];H)}^{2}+C\int_{0}^{t}\left\Vert \sigma (\tau )\right\Vert
_{H}^{2}d\tau \pier{,}
\end{eqnarray*}%
where $C_{01}:=C\left( \piergab{1+ {}}\left\Vert\pcol{\sigma _{0}^1}\right\Vert
_{H}^{2}+\left\Vert \piergab{s}_{1}\right\Vert^2_\infty\right)\!.$ Finally,
by using Gronwall lemma \pier{it is not difficult to conclude that}
\begin{eqnarray}
&&\left\Vert \sigma _{1}-\sigma _{2}\right\Vert _{C([0,t];H)\cap
L^{2}(0,t;V)}^{2}  \label{e19} \\
&\leq &C\left( \left\Vert \sigma _{0}^{1}-\sigma _{0}^{2}\right\Vert
_{H}^{2}+\left\Vert\piergab{s_{1}-s_{2}}\right\Vert _{L^{2}(0,t;H)}^{2}+\left\Vert
z_{1}-z_{2}\right\Vert _{C([0,t];H)}^{2}\right) ,\mbox{ for all }t\in
\lbrack 0,T],  \nonumber
\end{eqnarray}%
where $C$ stands for another constant depending additionally on the initial
data, $T$, and the source terms$.$

Next, we treat the initial-boundary value problem for the equation (\ref%
{e15-0}). To this end, \pier{for $z\in \mathcal{M}$ and $\sigma^z$ defined as above} 
we consider the intermediate equation%
\begin{equation}
\phi _{t}-\lambda \Delta \phi =g(\phi )f(z,\sigma ^{z},u),\mbox{ in }Q_{T},
\label{e20}
\end{equation}%
with 
\begin{equation}
\pier{\phi =0, \mbox{ on }\Sigma _{T};\ \ \mbox{ }\phi (0)=\phi _{0}, \mbox{ in }\Omega ,}
\label{e20-0}
\end{equation}%
\pier{where} 
\begin{equation}
g(r)=\left\{ 
\begin{array}{l}
-2r(1-r),\mbox{ }r\in [ 0,1], \\ 
0,\mbox{ \ \ \ \ \ \ \ \ \ \ \ \ \ otherwise.}%
\end{array}%
\right.  \label{e21}
\end{equation}%
We observe that $g$ is bounded and Lipschitz continuous on $\mathbb{R}$.
Recalling (\ref{f}) we note that $f(z,\sigma ^{z},u)\in L^{\infty }(Q_{T})$
and%
\begin{equation}
\left\Vert f(z,\sigma ^{z},u)\right\Vert _{\infty }\leq C(1+\left\Vert
u\right\Vert _{\infty }).  \label{e21-0}
\end{equation}%
We are going to prove that (\ref{e20})-(\ref{e20-0}) has a solution. We
apply the Banach fixed point theorem, \pier{by} fixing $\zeta \in C([0,T];H)$ and
studying the problem 
\begin{equation}
\phi _{t}-\lambda \Delta \phi =g(\zeta )f(z,\sigma ^{z},u),\mbox{ in }Q_{T},
\label{e21-1}
\end{equation}%
together with (\ref{e20-0}). We note that the right-hand side \pier{of (\ref{e21-1})} is in $%
L^{\infty }(Q_{T})$ and \pier{that} $\phi _{0}\in H.$ It is clear, \piergab{again by the 
solvability result for linear parabolic problems (see Ref.~ \cite{Lions1961}),} that 
\piergab{(\ref{e21-1}), (\ref{e20-0})} has a unique solution 
\[
\phi \in W^{1,2}(0,T;V_{0}^{\prime })\cap C([0,T];H)\cap L^{2}(0,T;V_{0}). 
\]

\pier{Now, we aim} at proving that the \pier{operator associating $\zeta$ to the solution $\phi $}
to (\ref{e21-1}) and (\ref{e20-0}) is a contraction \pier{mapping} on $C([0,T];H).$ For 
\pier{the sake of convenience}, we introduce the norm 
\[
\left\Vert \zeta \right\Vert _{B}:=\sup_{t\in \lbrack 0,T]}e^{-\gamma
t}\left\Vert \zeta (t)\right\Vert _{H},\mbox{ }\gamma >0,\mbox{ }\zeta \in
C([0,T];H), 
\]%
which is equivalent to the standard norm in $C([0,T];H).$

\pier{Let us} consider two solutions $\phi _{1}$ and $\phi _{2}$ to (\ref{e21-1}), (\ref%
{e20-0}) corresponding to $\zeta _{1}$ and $\zeta _{2},$ respectively. We
multiply the difference of \pier{the equations  (\ref{e21-1}) by $\phi _{1}-\phi _{2}$} and
integrate over $\Omega \times (0,t).$ With the help of \pier{(\ref{e20-0}) and (\ref{e21-0}) we have that} 
\begin{eqnarray*}
&&\frac{1}{2}\left\Vert (\phi _{1}-\phi _{2})(t)\right\Vert _{H}^{2}+\lambda
\int_{0}^{t}\left\Vert \nabla (\phi _{1}-\phi _{2})(\tau )\right\Vert
_{H}^{2}d\tau \\
&\leq &\int_{0}^{t}\int_{\Omega }\left\vert g(\zeta _{1})-g(\zeta
_{2})\right\vert \left\vert f(z,\sigma ^{z},u)\right\vert \left\vert \phi
_{1}-\phi _{2}\right\vert dxd\tau \\
&\leq &C\int_{0}^{t}\left\Vert (\zeta _{1}-\zeta _{2})(\tau )\right\Vert
_{H}\left\Vert (\phi _{1}-\phi _{2})(\tau )\right\Vert _{H}d\tau .
\end{eqnarray*}%
\pier{Hence, the} Young inequality and the Gronwall lemma allow us to get 
\[
\left\Vert (\phi _{1}-\phi _{2})(t)\right\Vert _{H}^{2}\leq
C\int_{0}^{t}\left\Vert (\zeta _{1}-\zeta _{2})(\tau )\right\Vert
_{H}^{2}d\tau . 
\]%
We multiply this inequality by $e^{-2\gamma t}$ and have successively%
\begin{eqnarray*}
&&e^{-2\gamma t}\left\Vert (\phi _{1}-\phi _{2})(t)\right\Vert _{H}^{2}\leq
Ce^{-2\gamma t}\int_{0}^{t}e^{2\gamma \tau }e^{-2\gamma \tau }\left\Vert
(\zeta _{1}-\zeta _{2})(\tau )\right\Vert _{H}^{2}d\tau \\
&\leq &Ce^{-2\gamma t}\int_{0}^{t}e^{2\gamma \tau }\left\Vert \zeta
_{1}-\zeta _{2}\right\Vert _{B}^{2}d\tau \leq \frac{C}{2\gamma }%
(1-e^{-2\gamma t})\left\Vert \zeta _{1}-\zeta _{2}\right\Vert _{B}^{2},\mbox{
for all }t\in \lbrack 0,T].
\end{eqnarray*}%
Taking the supremum with respect to $t\in \lbrack 0,T]$, we obtain%
\[
\left\Vert \phi _{1}-\phi _{2}\right\Vert _{B}\leq \sqrt{\frac{C}{2\gamma }}%
\left\Vert \zeta _{1}-\zeta _{2}\right\Vert _{B}
\]%
and note that for $\gamma >\frac{C}{2}$ the mapping \pier{$\zeta \mapsto \phi $}
is a contraction on $C([0,T];H)$. \pier{Thus,} we deduce that (\ref{e20})-(\ref%
{e20-0}) has a unique solution.

Next, we recall (\ref{e12-0}), that is $0\leq \phi _{0}\leq 1$ a.e. in $%
\Omega .$ Testing (\ref{e20}) by $-\phi ^{-}$ $(\phi ^{-}$ being the
negative part of $\phi )$, by a few \piergab{calculations} we obtain%
\begin{align*}
&\frac{1}{2}\left\Vert \phi ^{-}(t)\right\Vert _{H}^{2}+\lambda
\int_{0}^{t}\left\Vert \nabla \phi ^{-}(\tau )\right\Vert _{H}^{2}d\tau \\
&\leq 
\frac{1}{2}\left\Vert \phi _{0}^{-}\right\Vert
_{H}^{2}\pier{-}\int_{0}^{t}\int_{\Omega }\phi ^{-}g(\phi )f(z,\sigma
^{z},u)dxd\tau =0, 
\end{align*}
since \pier{$\phi_0 ^{-}=0$ a.e. in $\Omega$},
$\phi ^{-}=0$ a.e. in the set where $\phi \in \lbrack 0,1]$ and $%
g(\phi )=0$ a.e. in the set where $\phi \notin \lbrack 0,1].$ This implies
that $\phi ^{-}(t)=0,$ that is $\phi (t)\geq 0$ \pier{a.e. in $\Omega$,} 
for all $t\in \lbrack 0,T].$

Next, we multiply (\ref{e20}) by $(\phi -1)^{+}$ and by similar calculations
we deduce that%
\[
\frac{1}{2}\left\Vert (\phi -1)^{+}(t)\right\Vert _{H}^{2}+\lambda
\int_{0}^{t}\left\Vert \nabla (\phi -1)^{+}(\tau )\right\Vert _{H}^{2}d\tau
=0. 
\]%
Here, we used the fact that $\phi =0$ a.e. on $\Sigma _{T},$ so that $(\phi
-1)^{+}=0$ \pier{a.e.} on $\Sigma _{T}.$ Thus, we have \pier{shown} that 
$\phi (t)\leq 1$ \pier{a.e. in $\Omega$,} for all $t\in \lbrack 0,T].$

By these two last results we actually proved that $0\leq $ $\phi \leq 1$
a.e. on $Q_{T},$ implying that $g(\phi )=-2\phi (1-\phi ).$ It turns out
that this solution $\phi $ actually solves equation (\ref{e15-0}). Moreover,
we have that $\phi $ satisfies \pier{(cf.~(\ref{e8}))}
\begin{equation}
\phi \in W^{1,2}(0,T;V_{0}^{\prime })\cap C([0,T];H)\cap
L^{2}(0,T;V_{0})\cap L^{\infty }(Q_{T}).  \label{e21-2}
\end{equation}

Since \pier{the problem (\ref{e15-0}), (\ref%
{e20-0})} may also have other solutions fulfilling (\ref{e21-2}),
we have to prove a uniqueness result. For that, let us consider two
solutions $\phi _{1}$ and $\phi _{2}$ \pier{corresponding to} the same data, and compute the
difference of the \pier{respective equations  (\ref{e15-0}) tested} by $\phi _{1}-\phi
_{2}. $ \pier{Owing to a straightforward calculation and (\ref{e21-2}), we deduce~that} 
\begin{eqnarray*}
&&\frac{1}{2}\left\Vert (\phi _{1}-\phi _{2})(t)\right\Vert _{H}^{2}+\lambda
\int_{0}^{t}\pier{\left\Vert \nabla (\phi _{1}-\phi _{2})(\tau )\right\Vert _{H}^{2}}d\tau
\\
&=&-2\int_{0}^{t}\int_{\Omega }(\phi _{1}(1-\phi _{1})-\phi _{2}(1-\phi
_{2}))f(z,\sigma ^{z},u)(\phi _{1}-\phi _{2})dxd\tau \\
&\leq &2\left\Vert f(z,\sigma ^{z},u)\right\Vert _{\infty
}\int_{0}^{t}\int_{\Omega }(\phi _{1}-\phi _{2})^{2}\left\vert 1-(\phi
_{1}+\phi _{2})\right\vert dxd\tau \\
&\leq &C\int_{0}^{t}\left\Vert (\phi _{1}-\phi _{2})(\tau )\right\Vert
_{H}^{2}d\tau ,\pier{\mbox{
for all }t\in \lbrack 0,T],}
\end{eqnarray*}%
which by Gronwall lemma leads to the \pier{desired uniqueness property}.

The solution $\phi $ to (\ref{e15-0}), (\ref{e20-0}) we have found is
actually what we have denoted by $\phi ^{z}=\Psi (z).$ We have already shown
that $\phi ^{z}\in C([0,T];H)$ and that \pier{$0 \leq \phi ^{z}\leq1$} a.e. in $%
Q_{T}$\pier{, whence it follows} that $\Psi (\mathcal{M})\subset \mathcal{M}.$ 
\pier{Thus, for the last step it suffices to check} that $\Psi $ is a contraction mapping on $%
C([0,T];H)$.

First, we shall deduce a general estimate.

\pier{Let us} consider two pairs of data 
\pcol{$(\phi _{0}^i,\sigma _{0}^i,u_{i},\piergab{s_{i}};z_{i})$, $i=1,2,$} and the corresponding
solutions $(\phi _{i}=\phi ^{z_{i}},\sigma _{i}=\sigma ^{z_{i}}),$ $i=1,2.$
Then, we \pier{take} the difference of the equations (\ref{e15-0}), 
\pier{add $ \lambda(\phi _{1}-\phi _{2})$ to both sides, then
multiply} by $%
\phi _{1}-\phi _{2}$ and integrate over \pier{$\Omega\times (0,t) .$ Here, let} us denote 
$\phi :=\phi _{1}-\phi _{2},$ $\phi _{0}:=\phi _{0}^{1}-\phi _{0}^{2},$ \pier{and take into account} that $f$ is Lipschitz continuous with
a Lipschitz constant $L_{f}$ depending on $M,$ $m_{ref},$ $\rho ,$ $A,$ $%
\sigma _{r}$. We have 
\begin{eqnarray*}
&&\frac{1}{2}\left\Vert \phi (t)\right\Vert _{H}^{2}+\lambda
\int_{0}^{t}\left\Vert \phi (\tau )\right\Vert _{V}^{2}d\tau \\
&\leq& \frac{1}{2}%
\left\Vert \phi _{0}\right\Vert _{H}^{2}+\lambda \int_{0}^{t}\left\Vert \phi
(\tau )\right\Vert _{H}^{2}d\tau+2\int_{0}^{t}\int_{\Omega } \pier{\phi^2 \left\vert 1-\phi _{1}-\phi
_{2}\right\vert \left\vert f(z_{1},\sigma
_{1},u_{1})\right\vert} dxd\tau \\
&&+2\int_{0}^{t}\int_{\Omega }\left\vert f(z_{1},\sigma
_{1},u_{1})-f(z_{2},\sigma _{2},u_{2})\right\vert \left\vert \phi
_{2}(1-\phi _{2})\right\vert\pier{ \left\vert \phi\right\vert}
dxd\tau \\
&\leq &\frac{1}{2}\left\Vert \phi _{0}\right\Vert _{H}^{2}+\pier{\left( \lambda
 +4\left\Vert
f(z_{1},\sigma _{1},u_{1})\right\Vert _{\infty}\right) \int_{0}^{t}\left\Vert \phi (\tau )\right\Vert _{H}^{2}d\tau }\\
&&+2L_{f}\int_{0}^{t}(\left\Vert (z_{1}-z_{2})(\tau )\right\Vert
_{H}+\left\Vert (\sigma _{1}-\sigma _{2})(\tau )\right\Vert _{H}+\left\Vert
(u_{1}-u_{2})(\tau )\right\Vert _{H})\left\Vert \phi (\tau )\right\Vert
_{H}d\tau .
\end{eqnarray*}%
Therefore, recalling (\ref{e21-0}) and using the Young inequality, we infer
that 
\begin{align*}
&\frac{1}{2}\left\Vert \phi (t)\right\Vert _{H}^{2}+\lambda
\int_{0}^{t}\left\Vert \phi (\tau )\right\Vert _{V}^{2}d\tau 
\nonumber \\
&{}\leq \frac{1}{2}\left\Vert \phi _{0}\right\Vert
_{H}^{2}+C\int_{0}^{t}\left\Vert \phi (\tau )\right\Vert _{H}^{2}d\tau 
+2L_{f}\int_{0}^{t}\left\Vert (z_{1}-z_{2})(\tau )\right\Vert
_{H}^{2}d\tau
\\
&\quad {}+2L_{f}\left( \int_{0}^{t}\left\Vert (\sigma _{1}-\sigma _{2})(\tau
)\right\Vert _{H}^{2}d\tau +\int_{0}^{t}\left\Vert (u_{1}-u_{2})(\tau
)\right\Vert _{H}^{2}d\tau \right) ,
\end{align*}%
where $C$ depends on $\left\Vert u_{1}\right\Vert _{\infty }$ as well. Next,
we use (\ref{e19}) and obtain%
\begin{eqnarray*}
&&\frac{1}{2}\left\Vert \phi (t)\right\Vert _{H}^{2}+\lambda
\int_{0}^{t}\left\Vert \phi (\tau )\right\Vert _{V}^{2}d\tau \\
&&{}\leq \frac{1}{2}\left\Vert \phi _{0}\right\Vert
_{H}^{2}+C\int_{0}^{t}\left\Vert \phi (\tau )\right\Vert _{H}^{2}d\tau
+C\int_{0}^{t}\left\Vert z_{1}-z_{2}\right\Vert _{C([0,\tau ];H)}^{2}d\tau \\
&&\quad{}+C\left( \left\Vert \sigma _{0}^{1}-\sigma _{0}^{2}\right\Vert
_{H}^{2}+\left\Vert\piergab{s_{1}-s_{2}}\right\Vert _{L^{2}(0,T;H)}^{2}+\left\Vert
u_{1}-u_{2}\right\Vert _{L^{2}(0,T;H)}^{2}\right) .
\end{eqnarray*}%
Using the Gronwall lemma and estimating the right-hand side, we deduce that
\pier{%
\begin{eqnarray*}
&&\left\Vert \phi (s)\right\Vert _{H}^{2}+\int_{0}^{s}\left\Vert \phi (\tau
)\right\Vert _{V}^{2}d\tau 
%\leq C\int_{0}^{s}\left\Vert z\right\Vert_{C([0,\tau ];H)}^{2}d\tau \\
%&&+C(\left\Vert \phi _{0}\right\Vert _{H}^{2}+\left\Vert \sigma
%_{0}^{1}-\sigma _{0}^{2}\right\Vert _{H}^{2}+\left\Vert
%S_{1}-S_{2}\right\Vert _{L^{2}(0,T;H)}^{2}+\left\Vert u_{1}-u_{2}\right\Vert
%_{L^{2}(0,T;H)}^{2}) 
\\
&&{}\leq C\int_{0}^{t}\left\Vert z_{1}-z_{2}\right\Vert _{C([0,\tau ];H)}^{2}d\tau
+C \left( \left\Vert \phi _{0}\right\Vert _{H}^{2}+\left\Vert \sigma
_{0}^{1}-\sigma _{0}^{2}\right\Vert _{H}^{2}\right) \\
&&\quad{}+C\left(\left\Vert \piergab{s_{1}-s_{2}}\right\Vert _{L^{2}(0,T;H)}^{2}+\left\Vert
u_{1}-u_{2}\right\Vert _{L^{2}(0,T;H)}^{2}\right),
\end{eqnarray*}%
}%
for all $0\leq s<t\leq T.$ Taking now the supremum with respect to $0\leq
s\leq t$, \pier{we have~that}
%and taking into account that the first term on the right-hand side
%is increasing with respect to $t,$ we get%
\begin{eqnarray}
&&\left\Vert \phi _{1}-\phi _{2}\right\Vert
_{C([0,t];H)}^{2}+\int_{0}^{t}\left\Vert (\phi _{1}-\phi _{2})(\tau
)\right\Vert _{V}^{2}d\tau \label{e22} \\
&&{}\leq C\int_{0}^{t}\left\Vert
z_{1}-z_{2}\right\Vert _{C([0,\tau ];H)}^{2}d\tau +C\left(\left\Vert \phi _{0}^{1}-\phi _{0}^{2}\right\Vert _{H}^{2}+\left\Vert
\sigma _{0}^{1}-\sigma _{0}^{2}\right\Vert _{H}^{2}\right) \nonumber \\
&&\quad {}+ C\left(\left\Vert
\piergab{s_{1}-s_{2}}\right\Vert _{L^{2}(0,T;H)}^{2}+\left\Vert u_{1}-u_{2}\right\Vert
_{L^{2}(0,T;H)}^{2} \right),  \qquad \nonumber
\end{eqnarray}%
for all $t\in \lbrack 0,T].$ \pier{At this point, if we choose}
the same initial data and the same \piergab{$u_i$ and $s_i$, $i=1,2$,} we obtain%
\begin{equation}
\label{e22-0}
\left\Vert \phi _{1}-\phi _{2}\right\Vert _{C([0,t];H)}^{2}\leq
C\int_{0}^{t}\left\Vert z_{1}-z_{2}\right\Vert _{C([0,\tau ];H)}^{2}d\tau . 
\end{equation}
\pier{Hence, in the light of the argument used before, by considering here the norm 
\[
|\!|\!| z |\!|\!|
=\sup_{t\in \lbrack 0,T]}e^{-\gamma t}\left\Vert z
\right\Vert _{C([0,t];H)},\mbox{ }\gamma >0,\mbox{ }z \in C([0,T];H), 
\]
and playing with (\ref{e22-0}), we plainly obtain 
\[
 |\!|\!|\phi _{1}-\phi _{2} |\!|\!| \leq \sqrt{\frac{C}{2\gamma }} |\!|\!|z_{1}-z_{2} |\!|\!| .
\]%
Then, if $\gamma >\frac{C}{2},$ the mapping $z\mapsto  \phi ^{z}=\Psi
(z)$} is a contraction on $\mathcal{M}$ and $\phi ^{z}=z$ is the unique fixed
point$.$ Thus, problem (\ref{e1})-(\ref{e2}) with the initial and boundary
conditions has a unique solution, with all properties inherited from the
intermediate solutions occurred in the proof.

\pier{Now, we can replace $z_{i}$ by $\phi _{i}$, $i=1,2$, in (\ref{e22}) and apply} the
Gronwall lemma \pier{to infer that} 
\begin{eqnarray*}
&&\left\Vert \phi _{1}-\phi _{2}\right\Vert _{C([0,T];H)\cap
L^{2}(0,T;V)}^{2} \\
&\leq& C\left(\left\Vert \phi _{0}^{1}-\phi _{0}^{2}\right\Vert
_{H}^{2}+\left\Vert \sigma _{0}^{1}-\sigma _{0}^{2}\right\Vert
_{H}^{2}
\piergab{{}+\left\Vert u_{1}-u_{2}\right\Vert _{L^{2}(0,T;H)}^{2}
+\left\Vert\piergab{s_{1}-s_{2}}\right\Vert _{L^{2}(0,T;H)}^{2}} \right) . \qquad
\end{eqnarray*}%
Combining this estimate with (\ref{e19}) \pier{finally leads to} 
\begin{eqnarray}
&&\left\Vert \phi _{1}-\phi _{2}\right\Vert _{C([0,T];H)\cap
L^{2}(0,T;V)}^{2}+\left\Vert \sigma _{1}-\sigma _{2}\right\Vert
_{C([0,T];H)\cap L^{2}(0,T;V)}^{2}  \label{e23}\\
&\leq &C\left(\left\Vert \phi _{0}^{1}-\phi _{0}^{2}\right\Vert
_{H}^{2}+\left\Vert \sigma _{0}^{1}-\sigma _{0}^{2}\right\Vert
_{H}^{2}
\piergab{{}+\left\Vert u_{1}-u_{2}\right\Vert _{L^{2}(0,T;H)}^{2}
+\left\Vert\piergab{s_{1}-s_{2}}\right\Vert _{L^{2}(0,T;H)}^{2}} \right).\qquad  \nonumber
\end{eqnarray}

Let us now pass to (\ref{e3}) in which the right-hand side is in $L^{\infty
}(Q_{T}).$ It is easily seen, relying on the arguments before for $\sigma $,
that \pier{the probem based on (\ref{e3}), with the initial datum $p_{0}\in H$ and 
the homogeneous Neumann boundary condition,} has a unique solution 
\[
p\in W^{1,2}(0,T;V^{\prime })\cap C([0,T];H)\cap L^{2}(0,T;V). 
\]

Moreover, for the difference of two solutions corresponding to different
data we have the estimate 
\begin{align}
&\left\Vert p_{1}-p_{2}\right\Vert _{C([0,T];H)\cap L^{2}(0,T;V)}^{2}
\label{e24} \\
&\leq \left\Vert \BBB p_{0}^1-p_{0}^2\EEE\right\Vert _{H}^{2}+C\left\Vert \phi
_{1}-\phi _{2}\right\Vert _{L^{2}(0,T;H)}^{2}  \nonumber \\
&\leq C\, \left(\left\Vert p _{0}^{1}-p _{0}^{2}\right\Vert
_{H}^{2}+\left\Vert \sigma _{0}^{1}-\sigma _{0}^{2}\right\Vert
_{H}^{2} +\left\Vert \phi _{0}^1-\phi _{0}^2\right\Vert_H^2 \right.
\nonumber \\
&\qquad\quad \left.\piergab{{}+\left\Vert u_{1}-u_{2}\right\Vert _{L^{2}(0,T;H)}^{2}
+\left\Vert\piergab{s_{1}-s_{2}}\right\Vert _{L^{2}(0,T;H)}^{2}} \right)
 \nonumber
\end{align}%
By \pier{combining} (\ref{e23}) and (\ref{e24}), we obtain the continuity of the
solution with respect to the data, as claimed by (\ref{e14}).

Next, we prove the boundedness of $\sigma $ and $p$ and the property (\ref%
{e13-2}). We rewrite \pier{eq.~(\ref{e2})~as} 
\[
\sigma _{t}-\eta \Delta \sigma = \piergab{ S_h(1-\phi) +(S_c -s) \phi} 
-\gamma _{h}\sigma -\pier{\gamma _{ch}}\sigma \phi \in L^{\infty }(0,T;H),\mbox{ } 
\]%
and consider $\sigma _{0}\in L^{\infty }(\Omega ).$ By \piergab{Theorem 7.1, p.~181, in 
Ref.~ \cite{Ladyzhenskaia1968}}, we find that 
\begin{equation}
\sigma \in L^{\infty }(Q_{T}).  \label{e18}
\end{equation}%
Let $\sigma _{0}(x)\geq 0,$ \piergab{$s(t,x)\leq S_c$} a.e. and multiply (\ref{e2}) by $%
-\sigma ^{-}.$ We have 
\begin{eqnarray*}
&&\frac{1}{2}\left\Vert \sigma ^{-}(t)\right\Vert _{H}^{2}+\eta
\int_{0}^{t}\left\Vert \nabla \sigma ^{-}(\tau )\right\Vert _{H}^{2}d\tau
+\gamma _{h}\int_{0}^{t}\left\Vert \sigma ^{-}(\tau )\right\Vert
_{H}^{2}d\tau \\
&\leq &\frac{1}{2}\left\Vert \sigma _{0}^{-}\right\Vert _{H}^{2}+\left\vert
\gamma _{ch}\right\vert \int_{0}^{t}\int_{\Omega }(\sigma ^{-})^{2}\phi
dxd\tau -\int_{0}^{t}\int_{\Omega }\piergab{( S_h (1-\phi) + (S_c - s) \phi)}\sigma ^{-}dxd\tau \qquad \\
&\leq &\left\vert \gamma _{ch}\right\vert \int_{0}^{t}\left\Vert \sigma
^{-}(\tau )\right\Vert _{H}^{2}d\tau ,
\end{eqnarray*}%
since $0\leq \phi \leq 1$ and $\piergab{(S_c - s)} \geq 0$ a.e. in $Q_{T}.$ By Gronwall lemma
we obtain that $\sigma ^{-}(t)=0,$ hence $\sigma (t)\geq 0$ a.e. in $\Omega,$ for all $t\in \lbrack 0,T].$

By a similar argument we deduce that $p\in L^{\infty }(Q_{T})$ and $p(t)\geq
0$ a.e. in $\Omega $, for all $t\in \lbrack 0,T].$

Finally, we prove the further regularity of the solution \pier{and 
the estimate~(\ref{e13}). The properties $\phi \in X_0$, $\sigma \in X$, $p \in X$
can be compared with  (\ref{e8}) and we can observe that, e.g., for $\sigma$ and $p$ 
the triplet of spaces $(H,V,W) $ now replaces $(V', H,V) $. This fits into the frame 
of the structural regularity of solutions to parabolic systems (see, e.g., \piergab{Ref.~ \cite{Lions1961}).}}

Let $\phi _{0}\in V_{0}$ and write (\ref{e1}) as 
\[
\phi _{t}-\lambda \Delta \phi =-2\phi (1-\phi )f(\phi ,\sigma ,u)\in
L^{2}(0,T;H). 
\]%
\pier{Since the} operator $B_{0}:W_{0}\subset H\rightarrow H,$ $B_{0}=-\pier{\lambda}\Delta $ is
potential (that is, $B_{0}\phi =\partial j(\phi ),$ $j:V_{0}\rightarrow \pier{[0,+\infty]},$ $%
j(\phi )=\frac{\pier{\lambda}}{2}\int_{\Omega }\left\vert \nabla \phi \right\vert ^{2}dx)$,
according e.g.~to \piergab{Theorem~4.18, p.~179, in Ref.~ \cite{Barbu2010},}  it
follows that 
\[
\phi \in W^{1,2}(0,T;H). 
\]%
\pier{Then, we can test (\ref{e1}) by $\phi _{t}$; taking into account that (cf.~(\ref{e1}))} 
$$
\left\Vert f(\phi ,\sigma ,u)\right\Vert _{\infty }\leq C(1+\left\Vert \phi
\right\Vert _{\infty }+\left\Vert u\right\Vert _{\infty }), $$
we have%
\[
\int_{0}^{t}\left\Vert \pier{\phi _t} (\tau )\right\Vert _{H}^{2}d\tau +\frac{%
\lambda }{2}\left\Vert \phi (\tau )\right\Vert _{V_{0}}^{2}\leq \frac{%
\lambda }{2}\left\Vert \phi _{0}\right\Vert
_{V_{0}}^{2}+C\int_{0}^{t}\pier{(1+\left\Vert 
u\right\Vert _{\infty })\left\Vert \phi _{t}(\tau)\right\Vert_{H}}d\tau , 
\]%
since $\left\Vert \phi \right\Vert _{\infty }\leq 1.$ This yields 
\begin{equation}
\left\Vert \phi _{t}\right\Vert _{L^{2}(0,T;H)}+\left\Vert \phi
(t)\right\Vert _{V_{0}}^{2}\leq C(\left\Vert \phi _{0}\right\Vert
_{V}^{2}+\left\Vert u\right\Vert _{L^{\infty }(Q_{T})}^{2}+1),\mbox{ for all 
}t\in \lbrack 0,T].  \label{e30}
\end{equation}%
By comparison in the equation for $\phi $ we \pier{realize that $\Delta \phi \in
L^2(0,T;H) $. Then, using the elliptic regularity theory, due to the homogeneous Dirichlet
boundary condition for $\phi$, we can conclude that} 
\begin{equation}
\left\Vert \phi \right\Vert _{L^{2}(0,T;W_{0})}^{2}\leq C(\left\Vert \phi
_{0}\right\Vert _{V}^{2}+\left\Vert u\right\Vert _{L^{2}(0,T;H)}^{2}+1).
\label{e31}
\end{equation}%
Thus, we have obtained $\phi \in X_{0}$.

Assume that $\sigma _{0}\in V$ and consider again the equation (\ref{e2})
written as 
\begin{equation}
\label{e30-0}
\sigma _{t}-\eta \Delta \sigma +\gamma _{h}\sigma =\piergab{S_h + ( S_{ch} - s) \phi} -\pier{\gamma _{ch}}
\sigma \phi \in L^{2}(0,T;H). 
\end{equation}
The operator \pier{$-\eta \Delta +\gamma _{h}I : W\subset H \to H$ ($I$ standing for the identity in $H$) is potential as well, and then $\sigma \in
W^{1,2}(0,T;H).$ Multiplying (\ref{e2}) by $%
\sigma _{t}$ and integrating over $\Omega \times (0,t)$, we easily arrive at}
\begin{equation}
\int_{0}^{t}\left\Vert \sigma _{t}(\tau )\right\Vert _{H}^{2}d\tau
+\left\Vert \nabla \sigma (t)\right\Vert _{H}^{2}\leq C(\left\Vert \sigma
_{0}\right\Vert _{V}^{2}+\left\Vert \piergab{s} \right\Vert _{L^{2}(0,T;H)}^{2}+1),%
  \label{e17}
\end{equation}%
for all $t\in \lbrack 0,T]$. Also, \pier{from a comparison in (\ref{e30-0}) it follows that $\Delta \sigma \in L^2(0,T;H) $. 
This property, together with the homogeneous Neumann boundary condition for $\sigma$, entails}
\begin{equation}
\left\Vert \sigma \right\Vert _{L^{2}(0,T;W)}^{2}\leq C(\left\Vert \sigma
_{0}\right\Vert _{V}^{2}+\left\Vert  \piergab{s}\right\Vert _{L^{2}(0,T;H)}^{2}+1).
\label{e17-1}
\end{equation}%
\pier{A similar argumentation for $p$ yields} 
\begin{equation}
\pier{\left\Vert p\right\Vert _{W^{1,2}(0,T;H) \cap C([0,T];V)  \cap L^{2}(0,T;W)}^{2}\leq C(\left\Vert
p_{0}\right\Vert _{V}^{2}+\left\Vert \phi\right\Vert _{L^{2}(0,T;H)}^{2}+1).}
\label{e33}
\end{equation}%
\pier{By collecting  (\ref{e30}), (\ref{e31}), (\ref{e17})-(\ref{e33}) it is easy to get (\ref%
{e13}),} as claimed. \pier{Let us add the comment that, under the additional regularity 
specified by the estimate (\ref{e13}), the triplet $(\phi,\sigma, p)$ turns out to be a strong solution
to the system (\ref{e2})-(\ref{e6}), that is, it actually satisfies equations and boundary conditions almost everywhere in $Q_T$ and on $\Sigma_T$, respectively. By this remark, we conclude our proof.}
\hfill $\square $

\section{Numerical method}

\subsection{Space discretization}\label{space_disc}

We begin by deriving a weak form of the governing equations, in preparation for an isogeometric discretization based on the Galerkin method. The weak form is derived multiplying the governing equations with suitable smooth functions, integrating over the domain, and performing integration by parts. 
%We define the space $\mathcal{V}=\mathcal{H}^1$, where $\mathcal{H}^1$ is the Sobolev space of square-integrable functions with square-integrable first derivatives. The space $\mathcal{V}_0$ is defined as $\mathcal{V}_0=\{ v\in\mathcal{V} \; | \; v=0 \text{ on } \partial\Omega \}$. 
We recall that the spaces $V=H^1(\Omega)$ and $V_0=H^1_0(\Omega)$, defined in Section~\ref{wellpos}.
The weak form can be stated as follows: find $\phi\in{V}_0$, $\sigma\in{V}$ and $p\in{V}$ such that for all $\{w_1,w_2,w_3\} \in V_0 \times V \times V$,
\begin{align} 
\label{wf1} & \int_{\Omega} w_1 \left[ \phi_t + 2\phi(1-\phi)f(\phi,\sigma,u) \right]dx + \int_{\Omega} \lambda\nabla w_1\cdot  \nabla \phi dx  = 0, \\
\label{wf2} & \int_{\Omega} w_2 \left[ \sigma_t + ( \gamma_h(1-\phi) + \gamma_c\phi ) \sigma - S_h(1-\phi) - (S_c-s)\phi \right] dx\\
&\nonumber \qquad {  } +  \int_{\Omega} \eta\nabla w_2\cdot  \nabla \sigma dx  = 0, \\
\label{wf3} & \int_{\Omega} w_3 \left[ p_t + \gamma_p p -\alpha_h (1-\phi) - \alpha_c\phi \right] dx + \int_{\Omega} D\nabla w_3\cdot  \nabla p dx  = 0.
\end{align}
Once the equations have been written in weak form, they can be discretized by restricting the weak form to finite dimensional spaces ${V}^h\subset{V}$ and ${V}_0^h\subset{V}_0$. In this work, we construct these spaces using quadratic B-spline basis functions \cite{hughes2005isogeometric}. As an example, we take ${V}^h=\hbox{\rm span}\{N_j\}_{j=1,\dots,n_b}$ where the $N_j$'s are multivariate B-splines. The function $\phi^h$, which is a finite dimensional approximation to $\phi$, is defined as $\phi^h(t,{x})=\sum_{j=1}^{n_b}\phi_j(t)N_j({x})$ where the $\phi_j$'s are referred to as control variables. The functions $\sigma^h$, $p^h$, $w_i^h$ for $i=1,2,3$ are defined analogously.
Because $\phi^h\in V^h_0$, the control variables $\phi_j$ corresponding to the control points $j$ on the boundary must be zero to satisfy the homogeneous Dirichlet boundary conditions.

\subsection{Time discretization}

We integrate in time using the generalized-$\alpha$ method \cite{chung1993time,jansen2000generalized}. To define our time integration scheme, we introduce the global vector of degrees of freedom $\boldsymbol{\Phi}=\{\phi_j\}_{j=1,\dots,n_b}$. Analogously, we define $\boldsymbol{\Sigma}$ and $\boldsymbol{P}$ as the global vectors of control variables for $\sigma^h$ and $p^h$, respectively. Let us introduce the following residual vectors
\begin{equation}
\boldsymbol{R}=\{ \boldsymbol{R}_j \}; \quad \boldsymbol{R}_j = \{R_{\phi,j}, R_{\sigma,j}, R_{p,j}  \},
\end{equation}
where
\begin{align} 
\label{r1} & R_{\phi,j} = \int_{\Omega} N_j \left[ \phi^h_t + 2\phi^h(1-\phi^h)f(\phi^h,\sigma^h,u) \right]dx + \int_{\Omega} \lambda\nabla N_j\cdot  \nabla \phi^h dx,   \\
\label{r2} & R_{\sigma,j} =\int_{\Omega} N_j \left[ \sigma^h_t + ( \gamma_h(1-\phi^h) + \gamma_c\phi^h ) \sigma^h  - S_h(1-\phi^h) - (S_c-s)\phi^h \right] dx \\
& \nonumber \qquad {}  +\int_{\Omega} \eta\nabla N_j\cdot  \nabla \sigma^h dx,   \\
\label{r3} & R_{p,j} =\int_{\Omega} N_j \left[ p^h_t + \gamma_p p^h -\alpha_h (1-\phi^h) - \alpha_c\phi^h \right] dx + \int_{\Omega}D \nabla N_j\cdot  \nabla p^h dx.
\end{align}
Let us call $\boldsymbol{U}_n=\{ \boldsymbol{\Phi}_n, \boldsymbol{\Sigma}_n, \boldsymbol{P}_n \}$ the time discrete approximation to the control variables at time $t_n$. Our time-integration algorithm can be defined as follows: given $\vec U_n$, $\dot{\vec U}_n$ and the current time step $\Delta t_n=t_{n+1}-t_n$, find $\vec U_{n+1}$, $\dot{\vec U}_{n+1}$ such that
\begin{eqnarray}
 & & \vec R(\dot{\vec U}_{n+\alpha_m},{\vec U}_{n+\alpha_f})=\vec{0},  \\
 & & \vec U_{n+1}=\vec U_n+\Delta t_n\dot{\vec U}_n+\gamma\Delta t_n(\dot{\vec U}_{n+1}-\dot{\vec U}_n),  \\
 & & \dot{\vec U}_{n+\alpha_m}=\dot{\vec U}_n+\alpha_m(\dot{\vec U}_{n+1}-\dot{\vec U}_n),  \\
 & & {\vec U}_{n+\alpha_f}={\vec U}_n+\alpha_f({\vec U}_{n+1}-{\vec U}_n). 
\end{eqnarray}
The parameters $\alpha_m$, $\alpha_f$ and $\gamma$ define the stability and accuracy of the algorithm. The generalized-$\alpha$ method can be proven to be second-order accurate and $A$-stable by taking $\rho_\infty\in[0,1]$ and
\begin{equation}
\alpha_m=\frac{1}{2}\left(\frac{3-\rho_\infty}{1+\rho_\infty}\right),\quad \alpha_f=\frac{1}{1+\rho_\infty},\quad \gamma=\frac{1}{2}+\alpha_m-\alpha_f.
\end{equation}
In our calculations, we take $\rho_\infty=1/2$. The resulting nonlinear system of equations is linearized using the Newton-Raphson algorithm. We converged all the individual residuals $\boldsymbol{R}_{\phi}$, $\boldsymbol{R}_{\sigma}$, and $\boldsymbol{R}_{p}$ up to a certain tolerance $\varepsilon_{NR}$. This strategy ensured that the fulfilment of the convergence criterion by the global residual $\boldsymbol{R}$ was not due to potentially different scaling of the equations of the system (\ref{r1})-(\ref{r3}).  Within each iteration of  the Newton-Raphson algorithm, we solved the corresponding linearized system by means of the GMRES \cite{Saad1986} algorithm with a diagonal preconditioner. We chose a certain tolerance $\varepsilon_{GMRES}$ and a maximum number of iterations as convergence criteria for the GMRES algorithm.

\section{Simulation study}
\subsection{Description and parameter selection}

We consider two different cases of PCa in this simulation study: a mild tumor and an aggressive tumor. These different tumor behaviors can be implemented by appropriately choosing $\rho$ and $A$ within the function $m(\sigma)$ in ($\ref{msigma}$). In general, tumor growth is usually characterized by a high proliferation index and a low apoptotis index \cite{Hanahan2011}, which would translate into choosing $|\rho|>|A|$ within the context of our model. Then, we may represent mild tumors with lower values of these parameters for which the difference $|\rho|-|A|$ is also low, while more aggressive tumors can be modeled with higher values and such that the difference $|\rho|-|A|$ is also larger. For the simulations in this study, we calibrated $\rho$ and $A$ using the proliferation and apoptosis rates previously reported in the literature \cite{Berges1995,Schmid1993} in equations \eqref{rhodef}-\eqref{Adef}. We further computationally calibrated $m_{ref}$ in \eqref{gfunct}-\eqref{msigma} to ensure $|m(\sigma)-m_{ref}u|<1/3$ (see Section \ref{tumoreq}).

We begin our simulation study by analyzing mild and aggressive PCa growth in an untreated scenario. The metabolic profile of PCa is known to vary as it progresses towards a more malignant disease \cite{Hanahan2011,Trock2011}. Thus, we further consider four additional scenarios of untreated mild and aggressive PCa growth. We explore the effects of an effective and a poor nutrient supply within the tumor, i.e., high and low $S_c$ in \eqref{eq_nut}, respectively. Moreover, we analyze the effect of tumor metabolism on its dynamics by choosing a larger or a smaller value of $\gamma_c$ in \eqref{eq_nut}. 

After studying untreated tumor growth, we proceed to analyze the effects of cytotoxic and antiangiogenic therapy alone and combined for both the mild and the aggressive tumor. As introduced in Section \ref{chemotherapy}, we adopt the docetaxel-based cytotoxic chemotherapeutic protocol that is commonly used in the clinical management of advanced PCa. Hence, we consider 10 equal doses $d_c=$ 75 mg/m$^2$ of docetaxel delivered every three weeks \cite{Mottet2018,Kelly2012,Eisenberger2012}. We set $\tau_c=5$ days, according to pharmacodynamic studies of this drug \cite{Baker2004,Tije2005}. We choose bevacizumab as antiangiogenic drug and assume that the antiangiogenic therapy also consists of 10 equal doses $d_a=15$ mg/kg every three weeks \cite{Antonarakis2012,Kelly2012,picus2011phase}.  We set $\tau_a=30$ days \cite{Gordon2001,Lu2008,Ferrara2004}. We computationally calibrated the value of the constants $\beta_c$ and $\beta_a$ to produce a meaningful chemotherapeutic response as observed in clinical literature. 

Table \ref{parameters} summarizes the parameter selection discussed so far in this Section. The values of the other parameters in the model have been previously justified in the literature \cite{lorenzo2016tissue,lorenzo2017hierarchically,lorenzo2019computer,Xu2016}.

\begin{table}[!h]
\setlength{\tabcolsep}{0pt}
\caption{Parameter selection.} 
{\small
{\begin{tabular}{p{0.42\linewidth}p{0.1\linewidth}p{0.2\linewidth}p{0.28\linewidth}}
\toprule
 Parameter & Notation & Value & Simulation case\\ 
\midrule
\multicolumn{4}{@{}l}{\textbf{Tumor dynamics}} \\
Diffusivity of the tumor phase field           & $\lambda$  & 640 \si{\micro\metre}$^2$/day & All\\ 
Tumor mobility & $M$   & 2.5 1/day  & All\\
Net proliferation scaling factor               & $m_{ref}$  & 7.55$\cdot 10^{-2}$ 1/day & All\\
Scaling reference for proliferation rate               & $\bar{K}_\rho$  & 1.50$\cdot 10^{-2}$ 1/day & All\\

\multirow[t]{2}{*}{Proliferation rate}         & \multirow[t]{2}{*}{$K_\rho$} 
                                                 & 0.8$\cdot 10^{-2}$  1/day & Mild tumor\\
                                               & & 1.50$\cdot 10^{-2}$  1/day & Aggressive tumor\\
Scaling reference for apoptosis rate               & $\bar{K}_A$  & 2.10$\cdot 10^{-2}$ 1/day & All\\                                               
\multirow[t]{2}{*}{Apoptosis rate}         & \multirow[t]{2}{*}{$K_A$} & 0.7$\cdot 10^{-2}$  1/day & Mild tumor\\
& & 1.37 $\cdot 10^{-2}$  1/day & Aggressive tumor\\

\multicolumn{4}{@{}l}{\textbf{Cytotoxic chemotherapy}} \\
Mean lifetime of cytotoxic drug          & $\tau_c$  & 5 day & Cytotoxic chemotherapy and combined therapy\\ 
Cytotoxic drug effect          & $\beta_c$  & 1.59$\cdot 10^{-2}$ 1/(mg/m$^2$) & Cytotoxic  chemotherapy and combined therapy\\
Cytotoxic drug dose          & $d_c$  &  75 mg/m$^2$ & Cytotoxic  chemotherapy and combined therapy\\

\multicolumn{4}{@{}l}{\textbf{Nutrient dynamics }} \\
Nutrient diffusivity                           & $\eta$     & $6.4 \cdot 10^{4}$ \si{\micro\metre}$^2$/day & All \\
Nutrient supply in healthy tissue         & $S_h$      & $2$ g/L/day & All\\
\multirow[t]{3}{*}{Nutrient supply  in tumor tissue}  & \multirow[t]{3}{*}{$S_c$} 
                                                 & 2.75  g/L/day & Reference\\
                                               & & 3.125 g/L/day & Rich supply\\
                                               & & 2.375 g/L/day & Poor supply\\
Nutrient uptake by healthy tissue    & $\gamma_h$ & $2$ g/L/day & All\\
\multirow[t]{3}{*}{Nutrient uptake by tumor tissue}  & \multirow[t]{3}{*}{$\gamma_c$} 
                                                 & 17    g/L/day & Reference\\
                                               & & 18    g/L/day & Higher metabolism\\
                                               & & 16    g/L/day & Lower metabolism\\
                                               
\multicolumn{4}{@{}l}{\textbf{Antiangiogenic therapy}} \\
Mean lifetime of antiangiogenic drug          & $\tau_a$  & 30 day & Antiangiogenic and combined therapy\\ 
Antiangiogenic drug effect           & $\beta_a$  & 0.04 g/L/day/(mg/kg) & Antiangiogenic and combined therapy\\
Antiangiogenic drug dose          & $d_a$  & 15 mg/kg & Antiangiogenic and combined therapy\\

\multicolumn{4}{@{}l}{\textbf{Tissue PSA dynamics}} \\
Tissue PSA diffusivity                         & $D$        &  640 \si{\micro\metre}$^2$/day & All\\
Healthy tissue PSA production rate             & $\alpha_h$ & 1.712$\cdot 10^{-2}$ ng/mL/cc/day & All\\
Tumoral tissue PSA production rate             & $\alpha_c$ & $\alpha_c=15\alpha_h$ & All\\
Tissue PSA natural decay rate                  & $\gamma_p$ & 0.274 1/day & All\\

\bottomrule
\end{tabular}}}
\label{parameters}
\end{table}

\subsection{Computational setup}

All simulations were run on a 2D square domain with side length of $L_d=3000$ \si{\micro\metre} and 256 isogeometric elements per side. 
We chose a constant time step of $\Delta t=0.1$ days in all simulations.
The convergence of the Newton-Raphson method was set at tolerance $\varepsilon_{NR}=10^{-3}$, while for the GMRES algorithm was set at $\varepsilon_{GMRES}=10^{-3}$ or a maximum of 500 iterations.

The initial tumor phase field is approximated as an ellipsoidal tumor placed in the center of the domain and semiaxes $a=150$ \si{\micro\metre} and $b=200$ \si{\micro\metre} parallel to the sides of domain. To implement this initial condition, we $L^2$-projected the hyperbolic tangent function
\begin{equation}
\phi_0(x,y) = 0.5 - 0.5\tanh\left( 10\left(\sqrt{\frac{(x-L_d/2)^2}{a^2} + \frac{(y-L_d/2)^2}{b^2}} - 1\right) \right)
\end{equation}
over the quadratic B-spline space supporting our spatial discretization. This operation provides the control variables $\phi_{0,A}$ of the spline representation of the phase-field initial condition, i.e., $\phi^h_0(\boldsymbol{x},t)=\sum_{A=1}^{n_b}\phi_{0,A}(t)N_A(\boldsymbol{x})$ (see Section \ref{space_disc}).

The initial conditions for the nutrient and the tissue PSA are approximations based on $\phi_0$ given by 
\begin{equation}
\sigma_0 = c^0_\sigma + c^1_\sigma\phi_0
\end{equation}
and
\begin{equation}
p_0 = c^0_p + c^1_p\phi_0.
\end{equation}
The constants $c^0_\sigma$, $c^1_\sigma$, $c^0_p$, and $c^1_p$ are computationally calibrated \cite{lorenzo2017hierarchically}, such that $\sigma_0$ and $p_0$ reproduce a constant value of the nutrient and tissue PSA within the tumor and the host tissue. Hence, we choose $c^0_\sigma=1$ g/L, $c^1_\sigma=-0.8$ g/L, $c^0_p=0.0625$ ng/mL/cc, and $c^1_p=0.7975$ ng/mL/cc.

In the simulations with the chemotherapeutic treatments, we initially let the tumor grow until $t=60$ days, when we deliver the first drug dose. Then, we apply the remaining 9 cycles of chemotherapy included in all treatment strategies (i.e., cytotoxic, antiangiogenic, and combined) until $t=249$ days, when the last drug dose is delivered. Finally, we let the simulation proceed until completing 1 year, i.e., $t=365$ days.

\subsection{Simulation of untreated prostate cancer growth}

\begin{figure}[!t]
\centerline{\includegraphics[width=\linewidth]{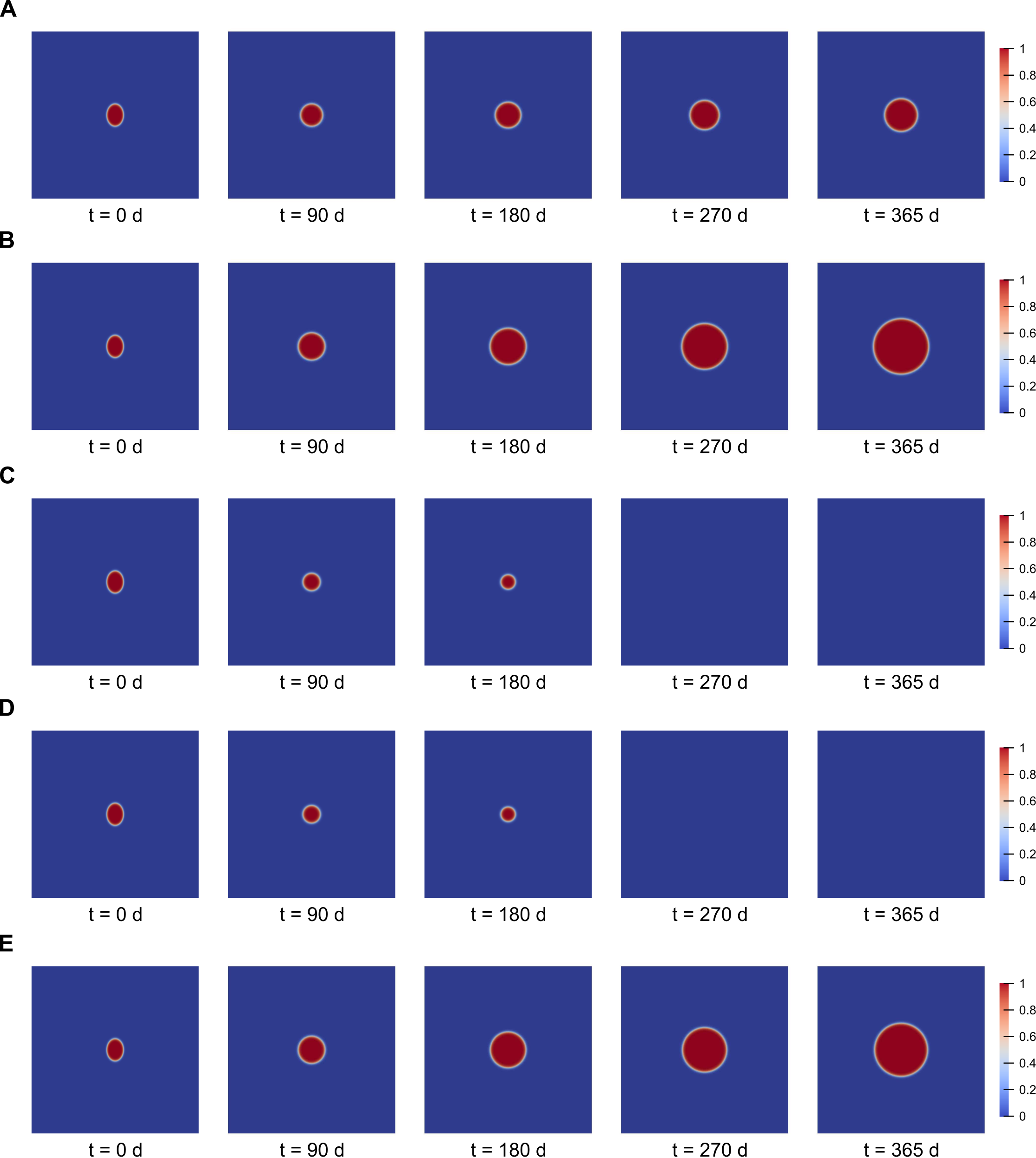}}
\vspace*{8pt}
\caption{Growth of a mild tumor under different assumptions of nutrient supply and tumor metabolism. (A) Reference simulation. (B) Larger nutrient supply within the tumor $S_c$. (C) Lower nutrient supply within the tumor $S_c$. (D) Larger tumor nutrient consumption rate $\gamma_c$. (E) Lower tumor nutrient consumption rate $\gamma_c$. }
\label{mild_tumor}
\end{figure}

\begin{figure}[!t]
\centerline{\includegraphics[width=\linewidth]{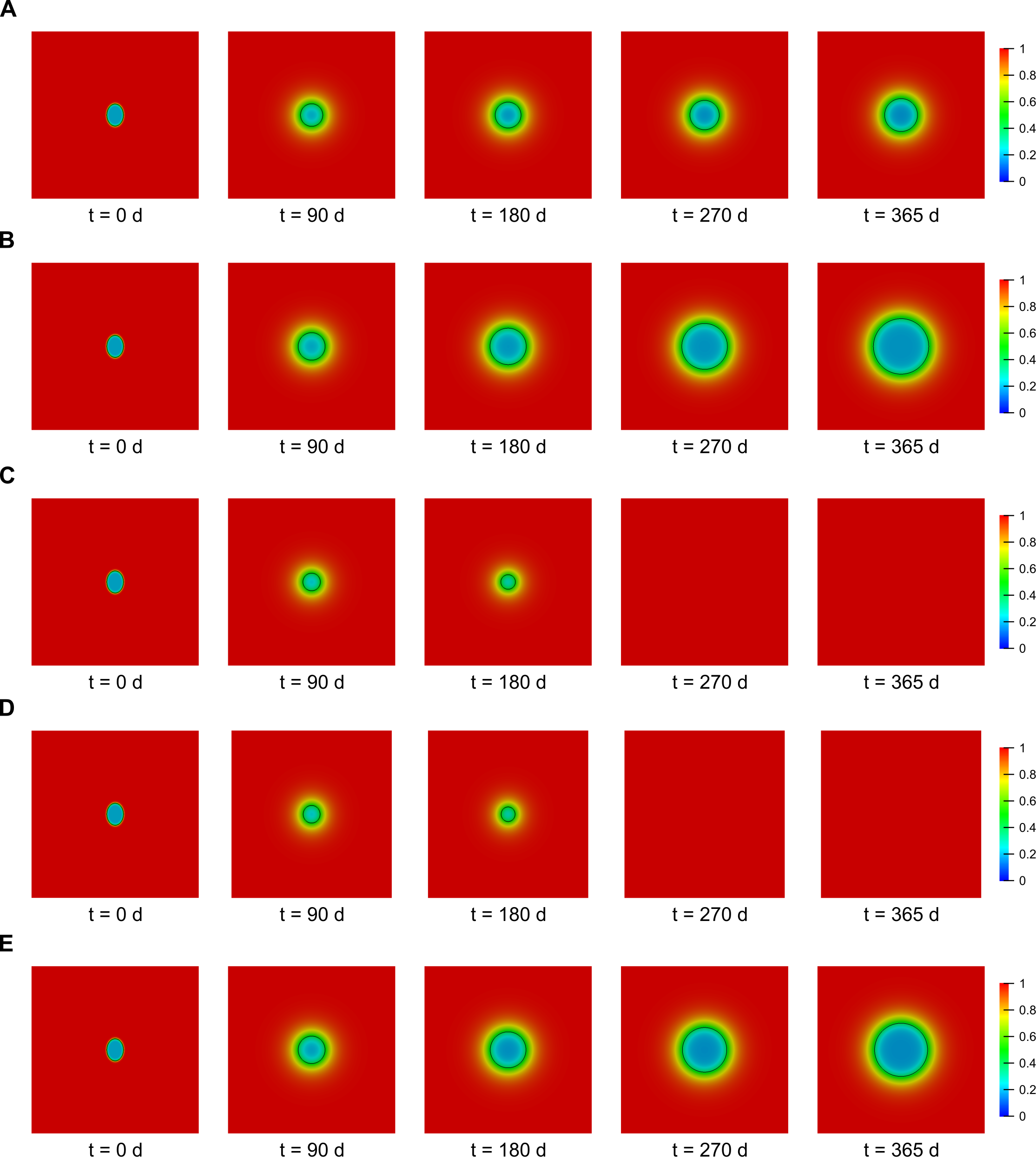}}
\vspace*{8pt}
\caption{Nutrient distribution during the growth of a mild tumor under different assumptions of nutrient supply and tumor metabolism. The tumor contour is depicted with a black line. (A) Reference simulation. (B) Larger nutrient supply within the tumor $S_c$. (C) Lower nutrient supply within the tumor $S_c$. (D) Larger tumor nutrient consumption rate $\gamma_c$. (E) Lower tumor nutrient consumption rate $\gamma_c$.}
\label{mild_nutrient}
\end{figure}

\begin{figure}[!t]
\centerline{\includegraphics[width=\linewidth]{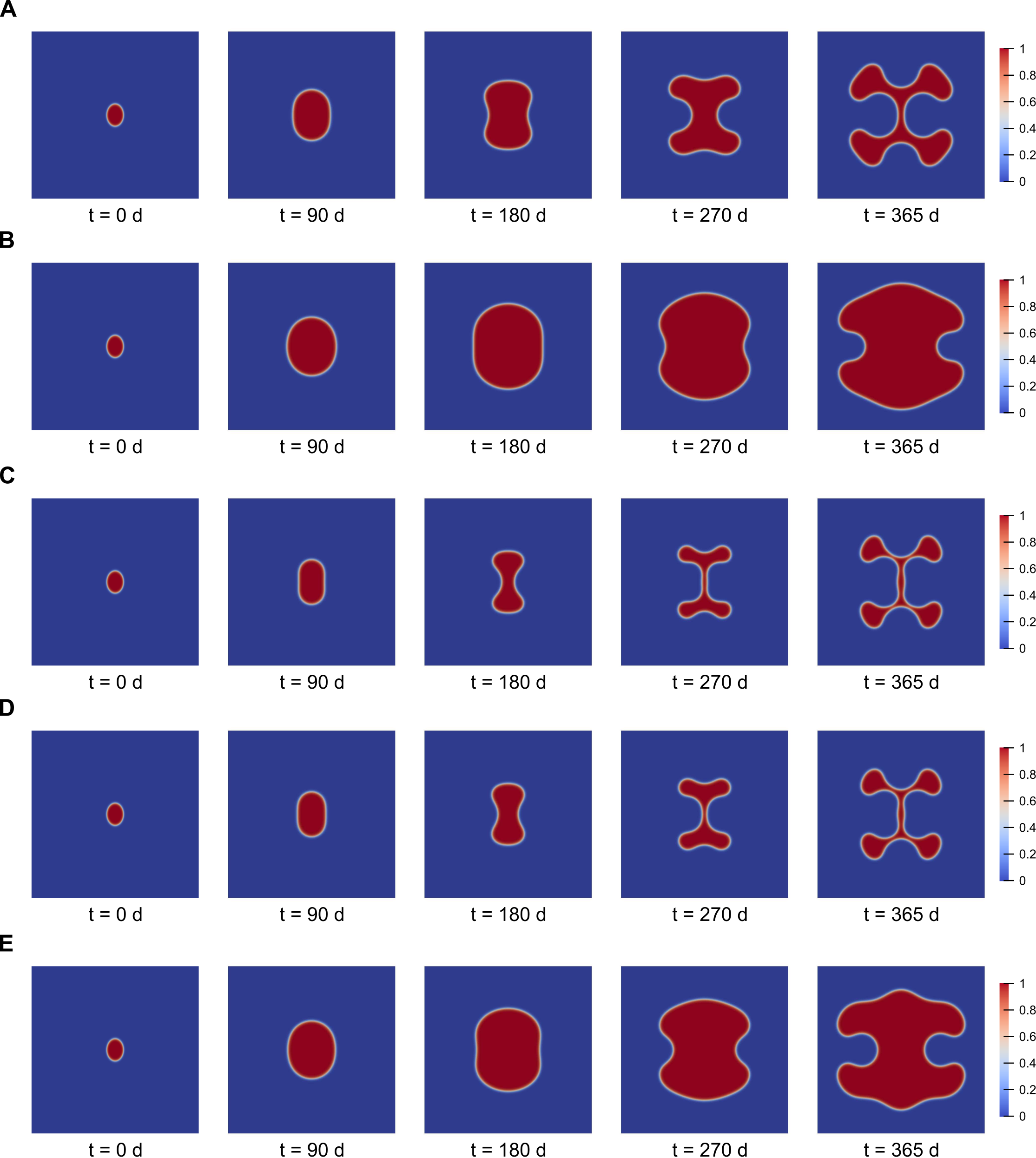}}
\vspace*{8pt}
\caption{Growth of an aggressive tumor under different assumptions of nutrient supply and tumor metabolism. (A) Reference simulation. (B) Larger nutrient supply within the tumor $S_c$. (C) Lower nutrient supply within the tumor $S_c$. (D) Larger tumor nutrient consumption rate $\gamma_c$. (E) Lower tumor nutrient consumption rate $\gamma_c$.}
\label{aggressive_tumor}
\end{figure}

\begin{figure}[!t]
\centerline{\includegraphics[width=\linewidth]{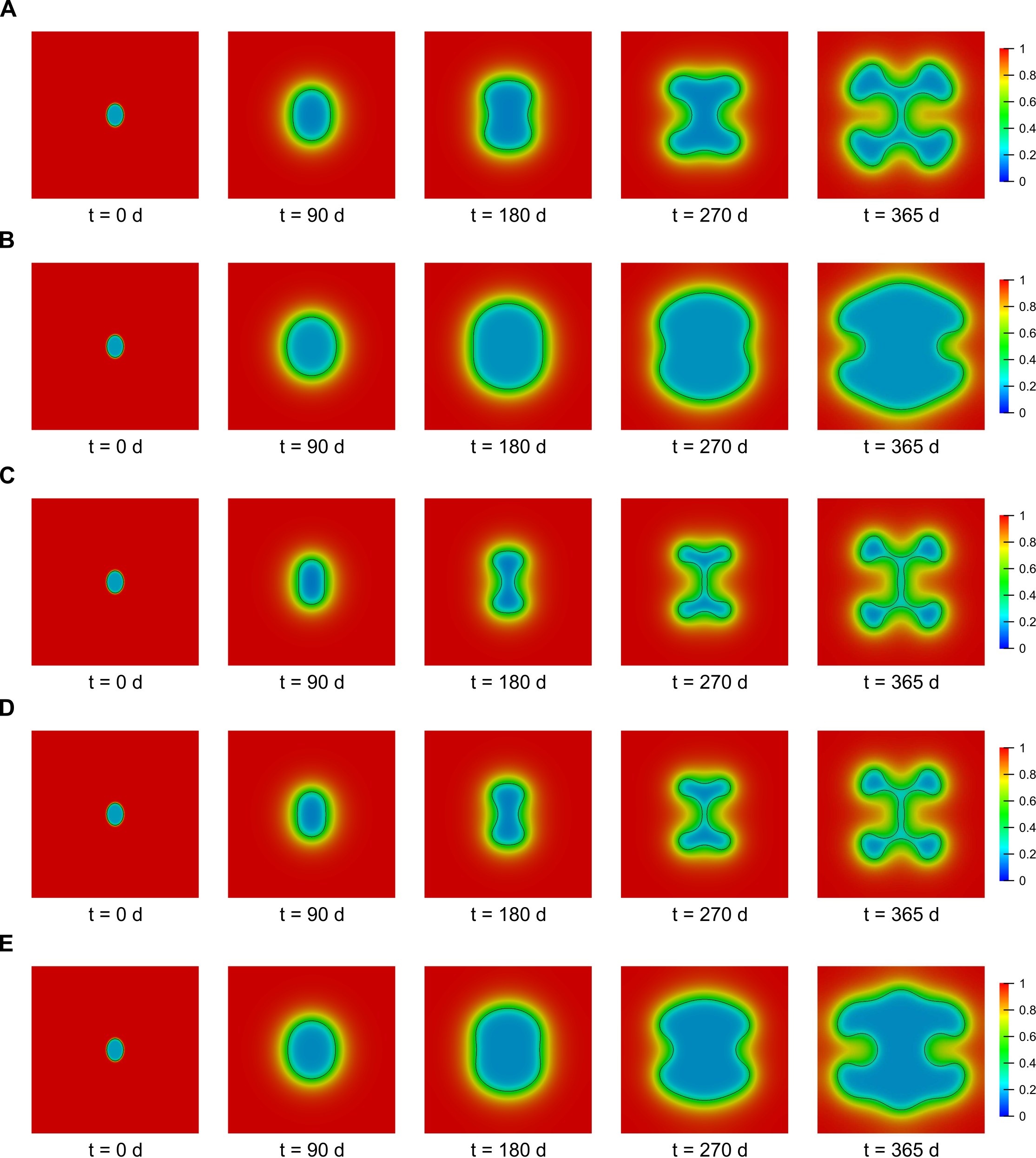}}
\vspace*{8pt}
\caption{Nutrient distribution during the growth of an aggressive tumor under different assumptions of nutrient supply and tumor metabolism. The tumor contour is depicted with a black line. (A) Reference simulation. (B) Larger nutrient supply within the tumor $S_c$. (C) Lower nutrient supply within the tumor $S_c$. (D) Larger tumor nutrient consumption rate $\gamma_c$. (E) Lower tumor nutrient consumption rate $\gamma_c$.}
\label{aggressive_nutrient}
\end{figure}

Figures \ref{mild_tumor} and \ref{mild_nutrient} respectively depict the evolution of the mild tumor and the corresponding nutrient distribution in all untreated simulation scenarios considered in this study.  We observe that the mild tumor grows at a slow pace and maintaining a round, spheroidal morphology that matches previous studies in the literature \cite{harma2010comprehensive,Erbersdobler2004,lorenzo2016tissue,Xu2016}. An increase in intratumoral nutrient supply or a metabolic adjustment decreasing the tumor nutrient uptake accelerates growth but does not alter the spheroidal morphology (see Figure \ref{mild_tumor}B,E). Conversely, a more scarce nutrient supply or a higher metabolic dependence on the nutrient arrests tumor growth and progressively lead to its extinction (see Figure \ref{mild_tumor}C,D). These results also align with the metabolic profile of mild prostatic tumors \cite{Hanahan2011,Trock2011}, which have not acquired a more competitive metabolism to support a highly proliferative behavior and still show a considerable dependence on the baseline nutrient in the host tissue. 

The growth of the aggressive tumor and the corresponding nutrient distribution in all untreated simulation scenarios are depicted in Figures \ref{aggressive_tumor} and \ref{aggressive_nutrient}, respectively. The aggressive tumor initially grows with a round morphology, but it eventually undergoes a morphological transformation by which it develops multiple branches and hence becomes more invasive \cite{harma2010comprehensive,Erbersdobler2004,noguchi2000assessment}. This phenomenon has also been observed in previous studies \cite{lorenzo2016tissue,Garcke2016,Wise2008}. Our simulations suggest that this shape instability is a tumor response to escape starvation because branching facilitates the access to the nutrient and leads to a nutrient redistribution that limits minimal concentrations in the inner regions of the tumor, which would potentially hamper PCa growth (see Figure \ref{aggressive_nutrient}). Hence, branching may be regarded as a mechanism of aggressive tumors to adapt to their harsh local environment.
The finger-like structures are thicker when we increase the intratumoral nutrient supply or reduce the metabolic dependence of the tumor on the considered nutrient (see Figure \ref{aggressive_tumor}B,E), whereas they become thinner in the opposite scenarios (see Figure \ref{aggressive_tumor}C,D). We have also observed that the morphological shift takes place sooner when nutrient availability is scarce or the tumor's metabolic dependence on the nutrient is larger. Conversely, the simulations with higher nutrient supply within the tumor and lower tumor nutrient consumption rate show a later appearance of branching. 

\begin{figure}[!t]
\centerline{\includegraphics[width=\linewidth]{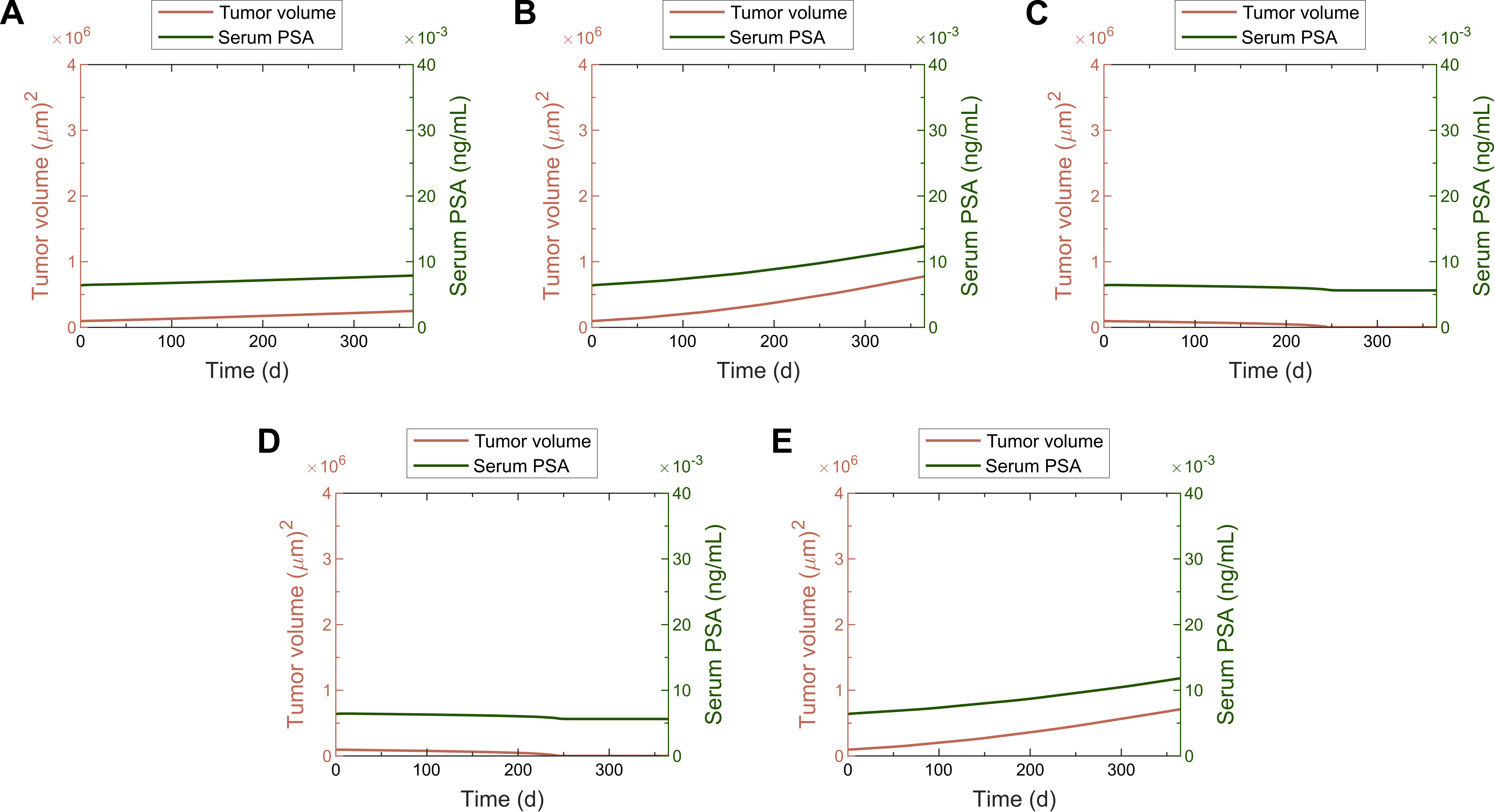}}
\vspace*{8pt}
\caption{Plots of tumor volume and serum PSA for the simulations with the untreated mild tumor under different assumptions of nutrient supply and tumor metabolism.(A) Reference simulation. (B) Larger nutrient supply within the tumor $S_c$. (C) Lower nutrient supply within the tumor $S_c$. (D) Larger tumor nutrient consumption rate $\gamma_c$. (E) Lower tumor nutrient consumption rate $\gamma_c$.}
\label{mild_uplots}
\end{figure}
\begin{figure}[!t]
\centerline{\includegraphics[width=\linewidth]{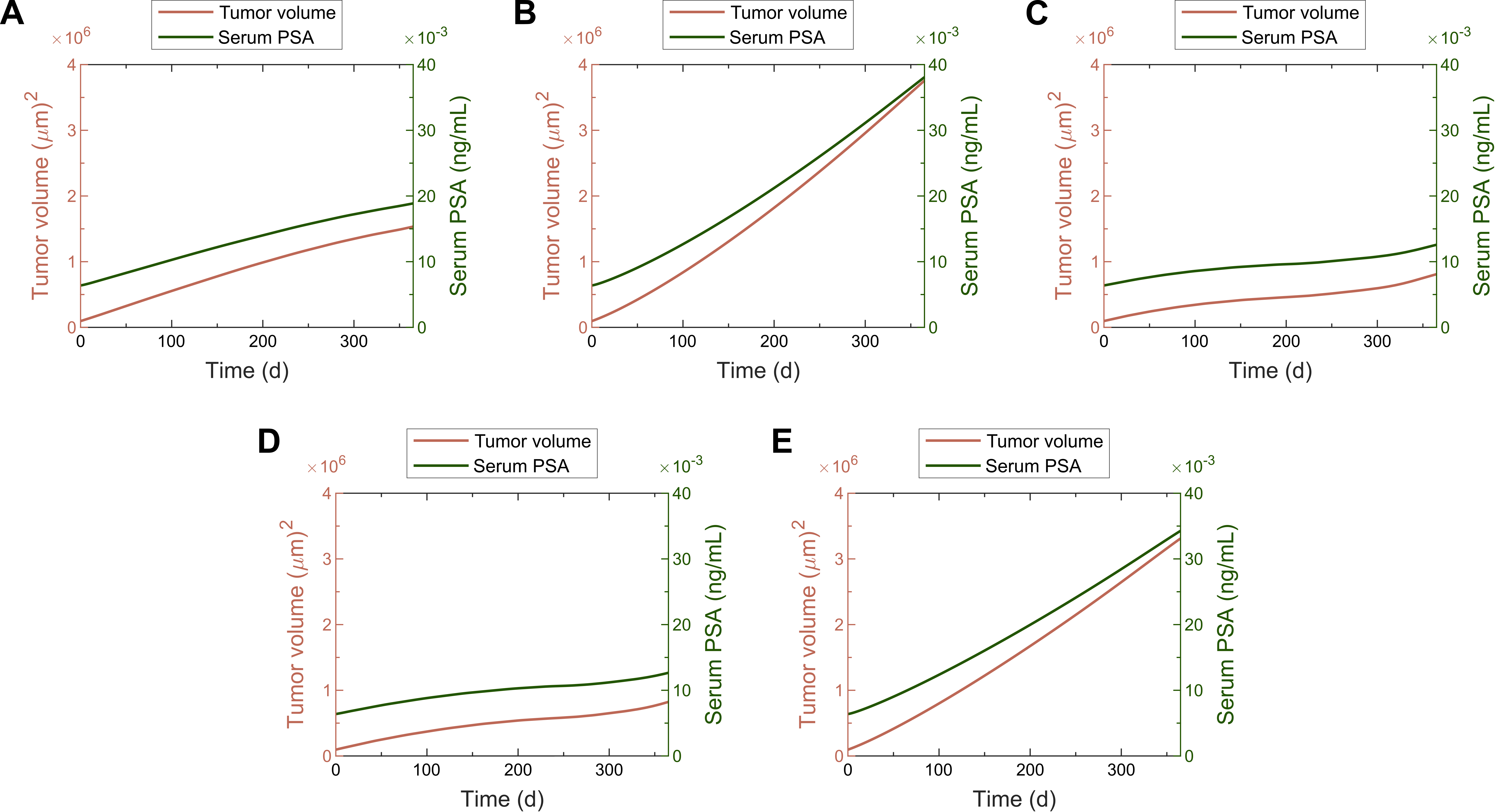}}
\vspace*{8pt}
\caption{Plots of tumor volume and serum PSA for the simulations with the untreated aggressive tumor under different assumptions of nutrient supply and tumor metabolism. (A) Reference simulation. (B) Larger nutrient supply within the tumor $S_c$. (C) Lower nutrient supply within the tumor $S_c$. (D) Larger tumor nutrient consumption rate $\gamma_c$. (E) Lower tumor nutrient consumption rate $\gamma_c$.}
\label{aggressive_uplot}
\end{figure}

Figure \ref{mild_uplots} displays the plots of tumor volume and serum PSA during each simulation of untreated mild PCa growth, while Figure \ref{aggressive_uplot} depicts the corresponding results for the aggressive cancer case.
In general, aggressive tumors grew faster and reached larger volumes in comparison to the corresponding simulation cases for the mild tumor, as shown in Figure \ref{aggressive_uplot}.
In the reference case of untreated aggressive PCa, the tumor volume increase slightly slowed down by the end (see Figure \ref{aggressive_uplot}A).
This trend was also observed in the first half of the simulations with decreased intratumoral nutrient supply and higher nutrient consumption by the tumor (see Figure \ref{aggressive_uplot}C,D). However, in these simulations we observe that the evolution of tumor volume exhibits a transient deceleration during branching. This marks an inflexion point that is later followed by a more rapid increase in tumor volume by the end of the simulation.
This faster dynamics corresponds to the thriving growth of the newly-formed branches with a more favorable nutrient distribution.
Hence, these observations suggest that branching was not completed in the reference simulation for the aggressive tumor, and the faster growing trend in tumor volume would ensue during the second year of simulation.

Figures \ref{mild_uplots} and \ref{aggressive_uplot} show that serum PSA followed a similar time evolution to tumor volume in all simulations of untreated PCa growth. This result is consistent with the extensive clinical use of this biomarker as a surrogate for the patient's tumor burden \cite{Mottet2018,lorenzo2016tissue,swanson2001quantitative,Vollmer2010}.
However, we observe that there were two cases in which tumor volume dynamics were noticeably faster than serum PSA dynamics, which correspond to the simulations featuring the aggressive tumor with increased intratumoral nutrient supply and lower nutrient uptake by cancerous tissue (see Figure \ref{aggressive_uplot}B,E).
Theses cases exhibited the fastest tumor growth dynamics as well as the largest tumor volumes registered in our simulation study and rapidly growing tumors have been previously reported to render comparatively lower serum PSA values \cite{swanson2001quantitative}.

\subsection{Simulation of treatment plans}

\begin{figure}[!t]
\centerline{\includegraphics[width=\linewidth]{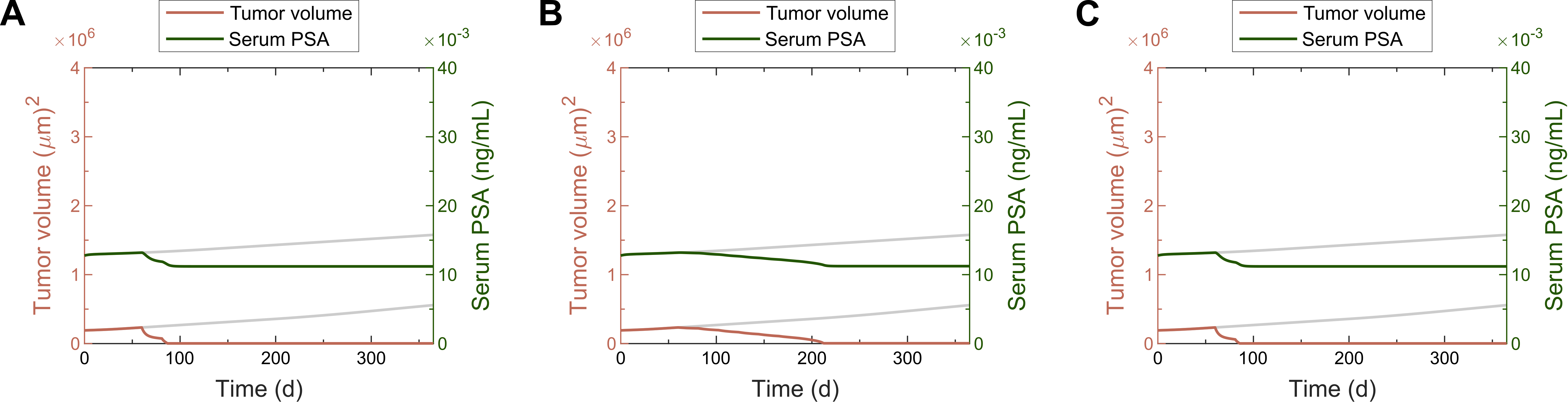}}
\vspace*{8pt}
\caption{Plots of tumor volume and serum PSA for the simulations with the mild tumor under different treatment plans. The gray lines in the background show the corresponding evolution of tumor volume and serum PSA in the untreated tumor reference scenario. (A) Cytotoxic chemotherapy. (B) Antiangiogenic therapy. (C) Combined therapy.}
\label{mild_ctplot}
\end{figure}

All treatment plans successfully removed mild PCa. Figure \ref{mild_ctplot} depicts the evolution of tumor volume and serum PSA under cytotoxic chemotherapy, antiangiogenic therapy, and combined therapy for the reference mild tumor. Cytotoxic chemotherapy was very effective, as it completely eliminated the tumor right after the second drug dose (see Figure \ref{mild_ctplot}A). Antiangiogenic therapy took comparatively longer to eliminate the tumor, as this required eight doses instead (see Figure \ref{mild_ctplot}B). We observe that the cytotoxic drug induces an immediate and steep decrease in tumor volume, which gradually smoothens between consecutive cytotoxic cycles as the concentration of the drug decays. Conversely, the antiangiogenic drug first slowly decelerates tumor growth and then drives its dynamics toward its extinction. Consequently, combined therapy was mostly driven by cytotoxic effects and, again, the tumor disappeared after the second cycle (see Figure \ref{mild_ctplot}C). Hence, our simulations show no advantage of combined therapy over cytotoxic chemotherapy alone for the mild tumor. 
%This outcome has also been observed in several clinical studies. 
In all simulations including drug effects, the reference mild tumor maintained a spheroidal, round morphology until vanishing and serum PSA followed similar dynamics to tumor volume evolution.

\begin{figure}[!t]
\centerline{\includegraphics[width=\linewidth]{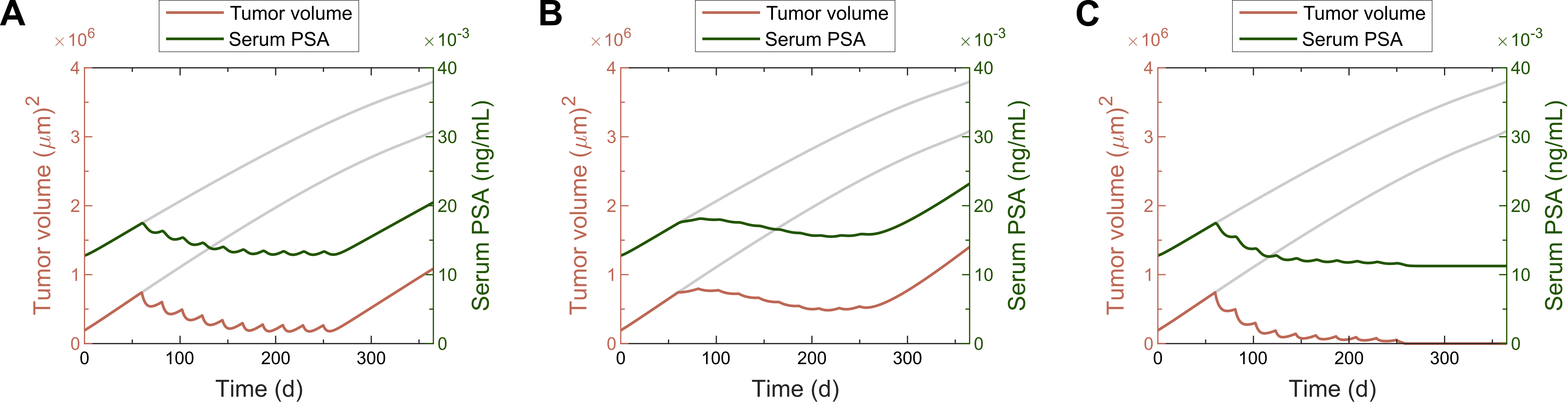}}
\vspace*{8pt}
\caption{Plots of tumor volume and serum PSA for the simulations with the aggressive tumor under different treatment plans. The gray lines in the background show the corresponding evolution of tumor volume and serum PSA in the untreated tumor reference scenario. (A) Cytotoxic chemotherapy. (B) Antiangiogenic therapy. (C) Combined therapy.}
\label{aggressive_ctplot}
\end{figure}

\begin{figure}[!t]
\centerline{\includegraphics[width=\linewidth]{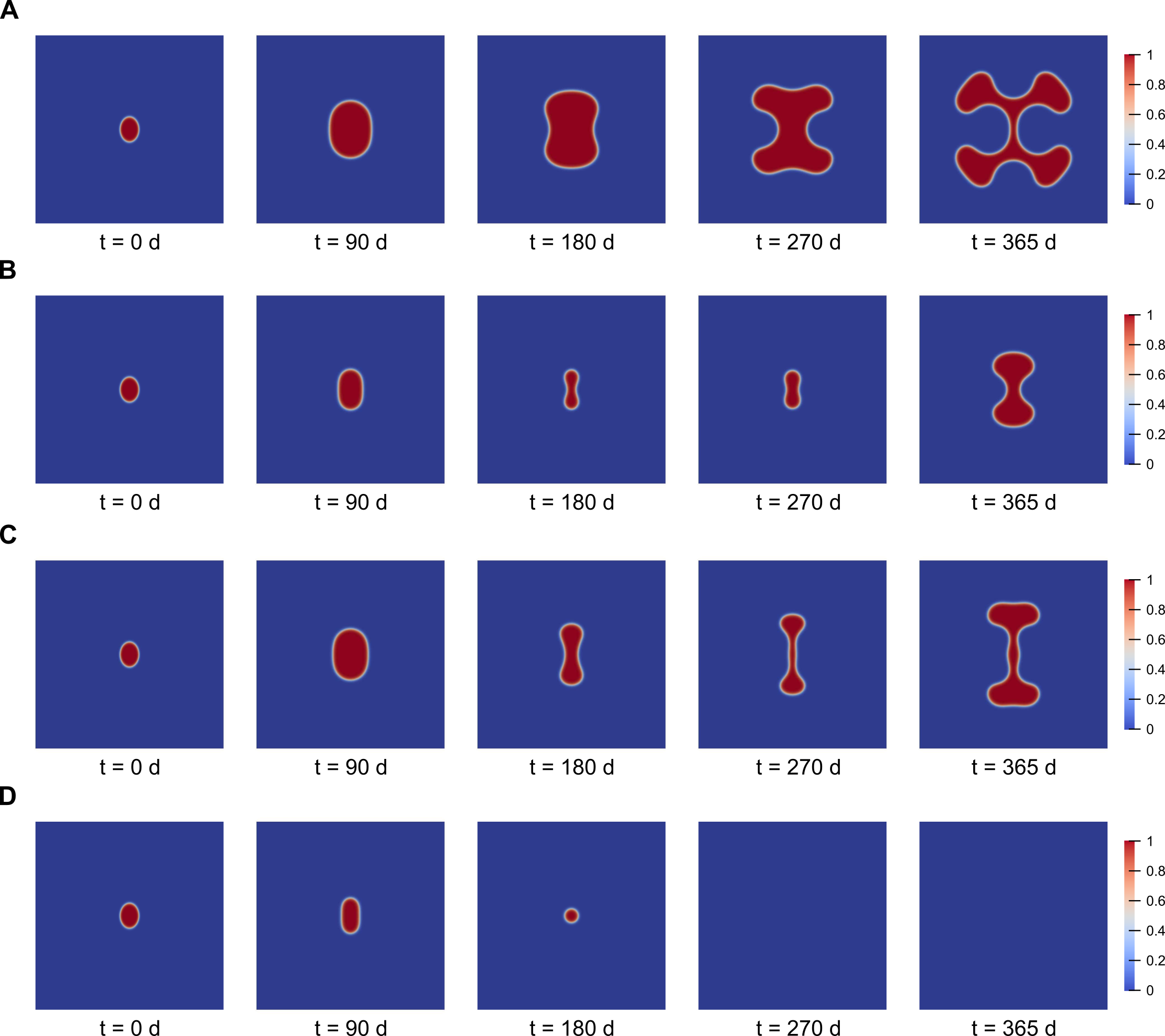}}
\vspace*{8pt}
\caption{Growth of an aggressive tumor under different treatment plans. (A) Untreated tumor. (B) Cytotoxic chemotherapy. (C) Antiangiogenic therapy. (D) Combined therapy.}
\label{ct_aggressive_tumor}
\end{figure}

\begin{figure}[!t]
\centerline{\includegraphics[width=\linewidth]{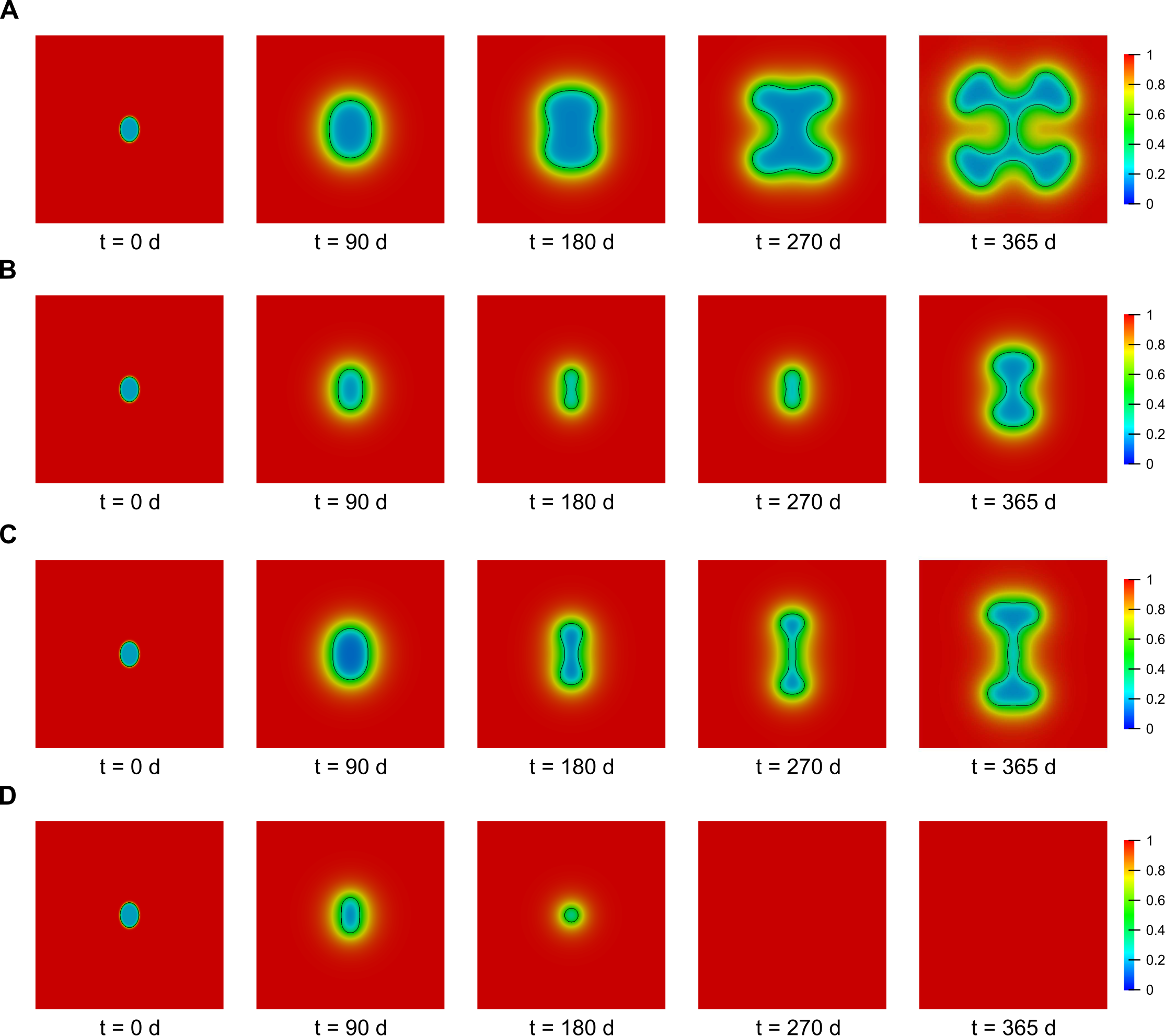}}
\vspace*{8pt}
\caption{Nutrient distribution during the growth of an aggressive tumor under different treatment plans. (A) Untreated tumor. (B) Cytotoxic chemotherapy. (C) Antiangiogenic therapy. (D) Combined therapy.}
\label{ct_aggressive_nutrient}
\end{figure}

Tumor volume and serum PSA dynamics during all simulated treatment plans for the reference aggressive tumor are plotted in Figure \ref{aggressive_ctplot}. Figures \ref{ct_aggressive_tumor} and \ref{ct_aggressive_nutrient} respectively show the morphological evolution of the aggressive tumor and the corresponding nutrient distribution during each treatment. Cytotoxic chemotherapy was able to temporarily control tumor growth and reduce the volume of the aggressive cancer, but did not completely eliminate it (see Figure \ref{aggressive_ctplot}A). Initially, each round of cytotoxic treatment produced an instantaneous and sheer decrease in tumor volume, as in the simulations with the mild tumor. However, the aggressive tumor managed to counteract this trend as the intensity of cytotoxic action decayed with the drug concentration and eventually resumed growth between consecutive cytotoxic cycles. This rendered a saw-like tumor volume evolution during the length of chemotherapy \cite{Kohandel2007,Powathil2007,Hahnfeldt2003,Henares-Molina2017,Kim2005,Caysa2012}. Serum PSA followed parallel dynamics to tumor volume. After the conclusion of the prescribed chemotherapeutic plan, the aggressive tumor resumed the growth dynamics exhibited before chemotherapy. 

Figure \ref{ct_aggressive_tumor}B shows that the tumor underwent branching during cytotoxic chemotherapy, leading to the development of two vertically connected cancerous masses by the end of the simulation. This result supports our interpretation of branching as a landmark feature of aggressive tumors. Figure \ref{ct_aggressive_nutrient}B suggests a possible mechanism that partially explains the failure of the simulated  cytotoxic chemotherapy. 
The small tumor volumes achieved during this treatment do not globally consume as much nutrient as larger untreated tumors. Hence, the nutrient concentration around and within the tumor is comparatively higher, which remarkably favors cancer growth once the cytotoxic drug concentration has decayed enough.
By the end of the chemotherapy plan, the cytotoxic drug is only able to reduce the tumor to a minimum limit volume, at which the tumor-promoting effect of the higher nutrient concentration overcomes the decaying cytotoxic action of the drug. This translates into the steady periodic pattern in tumor volume evolution and serum PSA dynamics during the last four cytotoxic cycles (see Figure \ref{aggressive_ctplot}A).

Antiangiogenic therapy was also unable to kill the aggressive cancer. While it provided a transient control on tumor volume too, it was less effective than cytotoxic chemotherapy (see Figure \ref{aggressive_ctplot}B). The first antiangiogenic drug dose slowed down tumor growth. Then, the ensuing ones induced a gradual and global decreasing trend in tumor volume, similar to that observed in the mild PCa case. The tumor-growth inhibitory effect of the antiangiogenic drug showed a little peak at the time of delivery of each dose, but it rapidly decreased afterwards and the tumor even grew momentarily between doses. Serum PSA exhibited similar dynamics to tumor volume. Figure \ref{ct_aggressive_tumor}C shows that the tumor developed two branches that grew vertically during the antiangiogenic treatment, without exhibiting further branching except until the end of the simulation.
This suggests that branching helped the tumor counteract the poorer intratumoral nutrient supply induced by antiangiogenic therapy, hence arising again as an adaptive feature of aggressive tumors.
Indeed, the morphology shift enabled the tumor to escape the inhibitory effect of antiangiogenic therapy before the end of treatment because tumor volume reaches a minimum between the eighth and ninth cycles and it exhibits a global increasing trend afterwards (see Figure \ref{aggressive_ctplot}B). Hence, antiangiogenic therapy had already failed before the delivery of the last two drug doses. Similar trends in tumor volume evolution during antiangiogenic therapy have been reported in the literature \cite{Benzekry2012,Hahnfeldt1999}.

Finally, combined therapy succeeded in eliminating the aggressive tumor right after the tenth cycle. Figure \ref{aggressive_ctplot}C shows that tumor volume evolution and serum PSA dynamics exhibited a global consistent decrease, even though both quantities briefly grew between consecutive doses when the aggressive tumor dynamics overcame the decaying combined tumor-inhibiting effect of the cytotoxic and antiangiogenic drugs.  The tumor initially evolved as an ellipsoid. During the third and fourth cycles, it quickly developed an ellongated vertical shape showing an incipient morphology shift. However, the tumor did not fully undergo branching and had reverted its morphology back to an ellipsoid by the fifth dose. During the second half of treatment and thereafter the tumor adopted a round, spheroidal morphology (see Figure \ref{ct_aggressive_tumor}D). We believe that the absence of branching in this simulation may be interpreted as an early sign of successful treatment outcome. Additionally, Figure \ref{ct_aggressive_nutrient}D shows that the nutrient concentration within small tumor volumes during combined therapy was lower than with cytotoxic chemotherapy. This is a consequence of the action of the antiangiogenic drug. Our simulations suggest that this reduction of intratumoral nutrient availability was a key enabling feature to kill the aggressive tumor.

\section{Conclusions}

Here, we present a mathematical model to describe the growth of prostatic tumors and their treatment using cytotoxic and antiangiogenic drugs. We describe the coupled dynamics of healthy and tumor tissue by leveraging the phase-field method and we assume that PCa growth is driven by a generic nutrient that follows reaction-diffusion dynamics \cite{gomez2018computational,lorenzo2016tissue,lorenzo2017hierarchically,lorenzo2019computer,Xu2016}. Our modeling approach enables the direct use of experimentally-determined tumor proliferation and apoptotis rates within the tumor phase-field equation. Cytotoxic chemotherapy is included as a tumor-inhibiting term in the phase-field equation. Antiangiogenic therapy is modeled as a term reducing the intratumoral nutrient supply in the nutrient dynamics equation. The adopted formulation for both treatment strategies aligns with previous mathematical models of chemotherapy and the pharmacodynamic behavior of common cytotoxic and antiangiogenic drugs used in PCa treatment.
To make the problem prostate-specific, we further couple the dynamics of the tumor to that of PSA \cite{lorenzo2016tissue,lorenzo2017hierarchically,swanson2001quantitative}, which is an extensively used biomarker of PCa in clinical management of the disease \cite{Mottet2018}.

\piergab{We prove that this model is well posed, that is, its solution exists and is unique in an appropriate functional space. In order to prove the existence of the solution, we apply a non-trivial argument based on the determination of fixed points in the systems of equations and initial and boundary conditions. Moreover, we derive further properties of the solution which are physically relevant. In particular, the solution component $\phi$, representing the phase variable, turns out to lie between $0$ (healthy phase) and $1$ (tumor phase), and the other components $\sigma$ and $p$ are non-negative and uniformly bounded, as expected from the model.}

Our simulations show that the model reproduces the common morphologies of untreated mild and aggressive tumors observed in previous computational, experimental, and clinical studies \cite{harma2010comprehensive,Erbersdobler2004,noguchi2000assessment,lorenzo2016tissue,lorenzo2019computer}.
We have also identified a branching instability as a feature characterizing a malignant behavior, which enables the tumor to adapt its morphology to facilitate the access to nutrients and avoid limit concentrations that would hamper its growth \cite{lorenzo2016tissue,Garcke2016,Wise2008}.
By varying the intratumoral nutrient supply and the tumor metabolic dependence on the nutrient, we were able to reproduce the starvation of a mild tumor and the adaptive branching morphologies of aggressive tumors driven by nutrient availability.
Both the evolution of tumor volume and serum PSA obtained in our simulations match the corresponding trends observed in clinical practice \cite{Mottet2018,Vollmer2010}.
Our model was also able to reproduce the poorer correlation between serum PSA and tumor volume in aggressive cancers exhibiting fast tumor dynamics \cite{swanson2001quantitative}.

Additionally, our results align with the use of cytotoxic chemotherapy as the reference approach to treat advanced PCa \cite{Mottet2018,Seruga2011,Eisenberger2012}, which may show further benefit from combination with an antiangiogenic drug in order to treat aggressive tumors \cite{Antonarakis2012,Seruga2011,Small2012}. Our simulations show that a mild tumor could be effectively treated with cytotoxic chemotherapy, antiangiogenic therapy, and combined therapy. However, our results suggest that cytotoxic drugs may suffice to effectively treat mild tumors.
Conversely, only combined therapy succeeded in eliminating the aggressive tumor in our simulations.
We observed that the combination of cytotoxic action with the reduction of intratumoral nutrient availability provided by antiangiogenic therapy was a pivotal mechanism to kill the aggressive tumor. Indeed, the tumor could not complete the morphology shift under combined therapy and reverted to a round, spheroidal tumor, which is more common in mild tumors. 
Cytotoxic and antiangiogenic therapy alone provided some tumor volume control, especially the cytotoxic approach. 
This reduction of tumor volume could be interesting in a neoadjuvant scenario, i.e., to optimize the outcome of a subsequent radical treatment (e.g., surgery, radiotherapy) \cite{Mottet2018}. 
Interestingly, our simulations suggest that the branching morphology shift may enable aggressive PCa to resist cytotoxic and antiangiogenic monotherapies. Chemoresistance is commonly attributed to a reduction in the drug effect on the tumor, but in our simulations we do not change the drug effect on PCa (i.e. $\beta_c$ and $\beta_a$). 
Hence, the study of tumor morphology may provide early information on the clinical outcome of the treatment with cytotoxic and antiangiogenic drugs. 

The mathematical model presented herein could be extended to include the evolving local angiogenic microvasculature \cite{figg2008angiogenesis,Hanahan2011}. 
This feature would improve the model by providing (1) a direct way to characterize intratumoral nutrient supply, (2) a direct target for antiangiogenic therapy, and (3) the possibility to refine the modeling action of cytotoxic and antiangiogenic therapy because these drugs would reach the tumor through this microvasculature network \cite{Kohandel2007,Powathil2007,Hahnfeldt1999,Hahnfeldt2003,Benzekry2013,Benzekry2012,Hinow2009}.
Hence, the introduction of tumor microvasculature would enable the study of drug supply to the tumor and would permit to explore the interactions between both forms of treatment considered herein, for instance, whether the reduction of microvasculature due to antiangiogenic therapy may significantly obstruct the effective delivery of cytotoxic drugs to the tumor and whether it would be more effective to normalize the tumor microvasculature to boost the supply of cytotoxic drugs \cite{Jain2005,Mpekris2017}.
Microvasculature may be modeled following alternative approaches \cite{Vilanova2017}. Hybrid models enable a precise description of the morphology and evolution of microvasculature, but they are more computationally expensive \cite{Frieboes2010,Vilanova2017,Xu2016}. Instead, a continuous reaction-diffusion equation describing the density of microvasculature would facilitate the computational coupling with our model and ensuing mathematical analysis at the cost of losing geometrical precision on the microvasculature \cite{Kohandel2007,Hahnfeldt1999,Benzekry2013,Benzekry2012,Vilanova2017}. 

Additionally, we plan to explore our mathematical model of PCa growth and cytotoxic and antiangiogenic drug effects within the context of optimal control problems \cite{Cavaterra2019,Colli2017,Garcke2018,Benzekry2013} to study the drug distributions that render an optimal treatment outcome. 
We also plan to explore \emph{in silico} the effects of cytotoxic and antiangiogenic drugs in 3D organ-scale, patient-specific scenarios using the same numerical approach described herein \cite{lorenzo2016tissue,lorenzo2017hierarchically,lorenzo2019computer}.
Computational models of tumor growth including the effect of mechanical deformation on cancer development have been found to provide superior predictions of pathological outcome \cite{Weis2015}. We have recently observed in a computational study that the mechanical stresses created by cancer growth and benign prostatic hyperplasia may obstruct prostatic tumor growth, which would explain the more favorable features of tumors arising in larger prostates \cite{lorenzo2019computer}. Thus, the extension of our model of PCa growth and drug effects to a poroelastic framework would not only improve the description of prostatic tumor dynamics, but it would also enable to include the effect of mechanical stress and fluid pressure on the delivery of nutrient and drugs to the tumor \cite{Fraldi2018,Jain2014,Roose2003}. This would permit to increase our understanding of the drug delivery and action in the complex tumor local environment as well as to refine current treatment strategies accordingly.
Our mathematical model could also be extended to a multiphase approach featuring various tumor species with varying sensitivities to the prescribed drugs \cite{Hahnfeldt2003,Jackson2000} and the modeling of the tumor-inhibiting effects of the drugs can be further refined \cite{Gorelik2008,Hinow2009}.
Radiation effects could be included in the tumor phase-field equation as an additional cytotoxic term \cite{Powathil2007,Corwin2013,Lima2017} within our modeling framework. 

Finally, we assumed that proliferation, apoptosis, and the drug effect rates $\beta_c$ and $\beta_a$ remained constant during the simulation. However, these parameters may evolve due to treatment action and the phenotypic evolution of the tumor \cite{Kim2005,Seruga2011,Hanahan2011}. Longitudinal series of PSA and medical images would enable the periodic update of these parameters during treatment, which could refine model predictions and provide early information on treatment outcome to guide further clinical decision-making (e.g., continuity of the treatment plan or shift to other treatment options). Additionally, this approach could contribute to extend our understanding of the mechanisms of chemoresistance in prostatic tumors and provide unique guidance in the design of optimal treatment protocols accordingly. 
In this context, our model of PCa growth with cytotoxic and antiangiogenic drug effects offers a mathematically robust framework with a vast modeling potential in order to explore personalized drug-based treatment strategies \emph{in silico}, which may assist physicians to successfully treat advanced PCa in the future.

\section*{Acknowledgements}
This research activity has been
performed in the framework of the Italian-Romanian collaboration agreement
\textquotedblleft Control and stabilization problems for phase field and
biological systems\textquotedblright\ between the Italian CNR and the
Romanian Academy.  The financial support of the project
Fondazione Cariplo-Regione Lombardia MEGAs-TAR \textquotedblleft Matematica
d'Eccellenza in biologia ed ingegneria come acceleratore di una nuova
strateGia per l'ATtRattivit\`{a} dell'ateneo pavese\textquotedblright\ is
gratefully acknowledged. The paper also benefits from the support of the
GNAMPA (Gruppo Nazionale per l'Analisi Matematica, la Probabilit\`{a} e le
loro Applicazioni) of INdAM (Istituto Nazionale di Alta Matematica) for PC
and ER. 
GL and AR have been also partially supported by the MIUR-PRIN project XFAST-SIMS (no. 20173C478N).
The present paper benefits from the support of the of
Ministry of Research and Innovation, CNCS --UEFISCDI, project number
PN-III-P4-ID-PCE-2016-0011, for GM.
The authors acknowledge the Rosen Center for Advanced Computing at Purdue University (USA) for providing HPC resources that contributed to the results presented in this paper.

\begin{footnotesize}
\bibliographystyle{abbrv}
%\section*{References}

\end{footnotesize}

\end{document}